\documentclass[a4paper,11pt,superscriptaddress]{article}

\usepackage[english]{babel}
\usepackage{amssymb}
\usepackage{amsfonts}
\usepackage{amsmath}

\usepackage{rotating}
\usepackage{tabulary} 
\usepackage{multirow}

\DeclareMathOperator\num{num}
\DeclareMathOperator\cnum{c-num}

\usepackage{graphicx}

\usepackage[round,comma,authoryear]{natbib}

\begin{document}

\title{On the limits of comparing subset sizes within $\mathbb{N}$}
\author{Sylvia Wenmackers\thanks{Centre for Logic and Philosophy of Science, Institute of Philosophy, KU Leuven, Belgium. E-mail: sylvia.wenmackers@kuleuven.be, URL: https://www.sylviawenmackers.be. \newline Forthcoming in the inaugural volume of the \textit{Journal for the Philosophy of Mathematics}.}}
\maketitle

\begin{abstract}
We review and compare five ways of assigning totally ordered sizes to subsets of the natural numbers: cardinality, infinite lottery logic with mirror cardinalities, natural density, generalised density, and $\alpha$-numerosity. Generalised densities and $\alpha$-numerosities lack uniqueness, which can be traced to intangibles: objects that can be proven to exist in ZFC while no explicit example of them can be given.
As a sixth and final formalism, we consider a recent proposal by \citet{Trlifajova:2024}, which we call c-numerosity. It is fully constructive and uniquely determined, but assigns merely partially ordered numerosity values. By relating all six formalisms to each other in terms of the underlying limit operations, we get a better sense of the intrinsic limitations in determining the sizes of subsets of $\mathbb{N}$.
\end{abstract}

\section{Introduction}\label{sec:intro}
The size of subsets of natural numbers played a central role in historical debates on the vexing notion of infinity.
Already in the ninth century, Thabit ibn Qurra defended the existence of different sizes of infinite collections of natural numbers \citep[see, e.g.,][§2]{Mancosu:2009}. \citet{Galileo:1638} famously compared the set of perfect squares to the full set of natural numbers. The finding that the former is both a proper subset of the latter (which suggests a smaller size) and can be put into one-to-one correspondence with it (which suggests equal size) is known as Galileo's paradox. Galileo concluded that sizes of infinite collections cannot be meaningfully compared.\footnote{On the first day of the dialogue, Salviati says that ``we cannot speak of infinite quantities as being the one greater or less than or equal to another'' \citep[p.~31]{Galileo:1638}; he argues for this based on a consideration of the `number' of all (natural) numbers and of (perfect) squares and their roots, from which he concludes that ``the attributes `equal,' `greater,' and `less,' are not applicable to infinite, but only to finite, quantities'' \citep[pp.~32--33]{Galileo:1638}.}

While many early mathematicians similarly rejected the idea that the infinite can be represented by an extended number system,\footnote{But see, e.g., \citet{Mancosu:2009} for notable exceptions, such as Grosseteste, Maignan, and Bolzano. We review a recent construction of Bolzano's approach in section~\ref{sec:c-num}.} later mathematicians have found different ways of expressing sizes of infinite sets. Here, we take `size' to be a pretheoretical term, that can be formalized in different ways: not only as cardinality, but also in terms of natural density and other approaches.
In this article, we take the set of natural numbers $\mathbb{N} \stackrel{\textrm{def}}{=} \{ 1, 2, 3, \ldots \}$ as the canonical example of a countably infinite set. We focus on its collection of subsets, the power set $\mathcal{P}(\mathbb{N})$, to discuss and compare different notions of size.

We start by reviewing and comparing five formalisms for assigning totally ordered (also called linearly ordered) sizes to subsets of $\mathbb{N}$: cardinality, infinite lottery logic, natural density, generalised density, and $\alpha$-numerosity.\footnote{These five formalisms suffice for our present purposes, without constituting an exhaustive overview. Other methods for assigning sizes to subsets of $\mathbb{N}$ exist. For instance, ideals and filters give a different notion of `small' and `large' sets \citep[see, e.g.,][p.~101]{Schechter:1997}. Although we will encounter ideals and filters in relation to $\alpha$-numerosities, and there may be interesting connections to other parts of this paper, we will not discuss this approach here.} In all cases, the notion of `size' is altered compared to its initial meaning that only applied to finite collections. This has been investigated before \citep[see, e.g.,][]{Mancosu:2009,Parker:2013}, but a new theory proposed by {\nobreak \citeauthor{Trlifajova:2024}} (\citeyear{Trlifajova:2024}) invites us to examine it once again. This paper highlights the non-uniqueness and non-constructive aspects of sizes assigned to infinite co-infinite sets on some accounts. In contrast, \citeauthor{Trlifajova:2024}'s (\citeyear{Trlifajova:2024}) theory is constructive but also requires us to drop yet another property of earlier notions of size.

The main goal of the present paper is not to defend a particular theory as providing us with the best notion of size, but rather to study the connections between various formalisms and to figure out which insights they give us collectively into size-related properties of subsets of $\mathbb{N}$. To achieve this, we will reconstruct the six theories in terms of an underlying limit operation that acts on the sequence of partial sums of the characteristic sequence of a given subset of $\mathbb{N}$ or the corresponding density, relative to $\mathbb{N}$.

We explore this issue in the context of standard set theory: Zermelo--Fraenkel set theory with the Axiom of Choice (ZFC), which takes the notion of arbitrary subsets of $\mathbb{N}$ (and other infinite sets) for granted.
Authors like \citet[p.~102]{Feferman:1999} objected to this approach, which raises a potential concern for our present project: does it really make sense to assign sizes to arbitrary subsets of natural numbers? Cardinality does so without a problem, which is unsurprising given that ZFC was designed with Cantor's approach in mind. Natural density fails to assign a size to some subsets of the set of natural numbers. Natural density can be extended to $\mathcal{P}(\mathbb{N})$, but not in unique way: the underdetermination shows up for some infinite co-infinite sets. As we will see, $\alpha$-numerosity is defined for all subsets of natural numbers but also lack uniqueness, which affects \emph{all} infinite co-infinite sets.

We discuss multiple ways of dealing with the underdetermination of $\alpha$-numerosity. In particular, \citeauthor{Trlifajova:2024}'s (\citeyear{Trlifajova:2024}) formalism trades the non-uniqueness of totally ordered $\alpha$-numerosity for a partially ordered notion of numerosity.

Throughout this article, we will pay special attention to \textit{intangibles}, in the terminology of \citet[§14.77]{Schechter:1997}: objects that can be proven to exist in ZFC but for which it can also be proven that no explicit example of them can be constructed, because their existence cannot be proven in ZF with a weaker choice principle, such as dependent choice (DC). \citet[§6.2 and §14.77]{Schechter:1997} describes ZF + DC as `quasiconstructive', in contrast to Bishop's constructivism, which rejects the law of the excluded middle. In this article, we similarly focus on the effects of invoking a nonconstructive choice axiom, without other constructivist considerations.\footnote{See \citet{Maschio:2020} for a discussion of natural density in the context of Bishop's constructivism.}
Like \citet[§14.77]{Schechter:1997}, we assume the view that intangibles are ``\emph{created by} our acceptance of the Axiom of Choice''.

Further on, it will be helpful to be able to compare how fine-grained different measures are. 
We say that a measure $m$ is \textit{more fine-grained} than a measure $M$ (or, equivalently, that $M$ is \textit{more coarse-grained} than $m$) if $m(S)=m(T)$ implies $M(S)=M(T)$ for all $S$ and $T$ in the domain where both measures are defined, but there exist $S$ and $T$ such that $M(S)=M(T)$ while $m(S) \neq m(T)$. In other words, there can be no $M$-difference without an $m$-difference, but the inverse does occur. This phrasing, familiar from the literature on emergence, already suggests that $M$ may be multiply realizable by various more fine-grained measures; this is indeed what we shall find.

The paper is structured as follows. Section~\ref{sec:sizes} reviews the five formalisms and compares them in terms of fine-grainedness and underdetermination. All formalisms are consistent on their own and have specific applications. Moreover, their results can be combined fruitfully, provided that we do not conflate their different notions of `size'. If we violate this condition, we run into inconsistencies such as Galileo's paradox. Section~\ref{sec:c-num} introduces the sixth approach by \citet{Trlifajova:2024} and reconstructs all theories in terms of different limit operations.
Section~\ref{sec:concl} summarizes our findings.

\section{Formalisms for assigning totally ordered sizes to subsets of $\mathbb{N}$}\label{sec:sizes}
\citet{Dedekind:1888} defined infinite sets as those sets that have a proper subset that can be put into one-to-one correspondence with the full set.\footnote{Bolzano had given a similar characterization earlier; see, e.g., \citet{Mancosu:2009}.} This shows that a version of Galileo's paradox will arise for any infinite set and thus complicates any assignment of size to such sets.
Set sizes can be compared by the subset relation (between a set and another one contained in it) and by one-to-one correspondence (between the elements of two sets of equal size). Both methods give compatible verdicts when comparing sizes of finite sets, but Galileo's paradox and Dedekind's definition both show that this is not the case for infinite sets.
This observation can be reconstructed as a dilemma between these two comparison methods \citep{Mancosu:2009}: 
to extend the notion of size to infinite sets, one can either respect the ordering induced by the subset relation (called the part--whole principle or the Euclidean principle\footnote{Named after the fifth and final `Common Notion' in Book~I of Euclid's \textit{Elements}, which says that the whole is greater than the part.}) or define equality based on the existence of a bijection (called Hume's principle\footnote{Named after \citet[Book~I, Part III, §I]{Hume:1739}: ``When two numbers [i.e., collections] are so combined, as that the one has always a unit answering to every unit of the other, we pronounce them equal'', as quoted (in German) by \citet[Ch.~IV, §63]{Frege:1884}.}), not both.

Historically, the second horn of the dilemma was chosen first: \citet{Cantor:1895} used it to develop cardinality theory. More recently, \citet{BenciDiNasso:2003b} chose the first horn to develop $\alpha$-numerosity theory. However, Bolzano's \textit{Paradoxes of the Infinite} from 1848 can be regarded as a precursor to the latter idea: \citet{Trlifajova:2024} has given a reconstruction, which we review in section~\ref{sec:c-num}.

As we will see, natural density retains neither: in number theory, it is usually interpreted as a probability \citep[see, e.g.,][Ch.~III.1]{Tenenbaum:2015}. If the measure does not express set size, then this classical dilemma does not apply to it. In this paper, however, we do consider natural density as a notion of (relative) size, showing that the paradox can be reconstructed as a trilemma instead.

Below, we review five formalisms that allow us to assign totally ordered sizes to subsets of $\mathbb{N}$: cardinality (§\ref{sec:card}), infinite lottery logic which includes mirror cardinalities (§\ref{sec:ILL}), natural density and generalised density (§\ref{sec:dens}), and $\alpha$-numerosity (§\ref{sec:num}).

\subsection{Cardinality}\label{sec:card}
Cardinality theory starts from Hume's principle, which takes the existence of a bijection between sets as indicative of their equal `size' expressed as a cardinal number. So, cardinality is defined via equivalence classes on sets between which there exists a one-to-one mapping.\footnote{For a brief contemporary presentation of cardinality theory, see, e.g., \citet[Ch.~3]{Jech:2003}.}

The cardinality of a (well-orderable) set $S$ is written as $|S|$. 
Assuming ZFC, all sets are well-orderable, so all sets have a well-defined cardinality. In particular, 
since $\mathbb{N}$ is well-ordered by its standard canonical order, cardinal numbers of subsets of $\mathbb{N}$ are fully determined: they are equal to the natural number of elements for finite subsets and they are equal to the first infinite (or `transfinite,' on Cantor's terminology) cardinal number, $\aleph_0$, for all infinite subsets of $\mathbb{N}$. In other words, while they are sensitive to singleton differences between finite sets (i.e., maximally fine-grained), they maximally coarse-grain the sizes of infinite sets by mapping co-finite as well co-infinite infinite sets to the same value.

\subsection{Mirror cardinality}\label{sec:ILL}
Rather than applying Hume's principle to sets themselves, one can apply it to their complements instead; doing so requires a fixed superset, here $\mathbb{N}$. This approach yields ``mirror images of the cardinals'' for co-finite sets, as discussed by \citet[p.~386]{Mancosu:2015}.

A related proposal has been raised in the context of probability. \citet{Norton:2021} considered full label-invariance on countable sets as a notion of uniformity stronger than mere singleton-uniformity: he required that the measure assigned to a subset of a countable set should be invariant under permutations of the labels. Since permutations are one-to-one mappings, we can expect the proposal to be close to cardinality.

Indeed, the approach of \citet{Norton:2021} distinguishes the probability of finite subsets by their finite cardinality (maximally fine-grained) and similarly distinguishes the probability of co-finite subsets. All infinite co-infinite sets get the same probability rank, so this remains coarse-grained. In particular, \citet[§8]{Norton:2021} introduced infinite lottery logic in terms of a chance function, $Ch$. To finite non-empty sets with $n$ elements, $Ch$ assigns the valuation $V_n$, which may be read as `unlikely'. To co-finite sets that are complementary to sets with $n$ elements, $Ch$ assigns the valuation $V_{-n}$, which may be read as `likely'. The values of the empty set and the full set are $V_0$, read as `certain not to happen' and $V_{-0}$, `certain to happen', respectively. Finally, $Ch$ maps all infinite co-infinite sets to the same value, $V_\infty$: this is read as `as likely as not'. The readings of the values are informal interpretations, motivated by their ordering. This is given by an antisymmetric, transitive and irreflexive relation, $<$, on them \citep[p.~S3866]{Norton:2021}:
$$V_0 < V_1 < V_2 < V_3 < \ldots < V_\infty < \ldots < V_{-3} < V_{-2} < V_{-1} < V_{-0}.$$

Here, we observe that this approach could be considered as an alternative conception of size as well, not based on the subset relation, but on the existence of bijections augmented with the relation of complements (relative to $\mathbb{N}$). Like cardinality, it is total on $\mathcal{P}(\mathbb{N})$ and uniquely determined, but it is more fine-grained. It still coarse-grains all infinite co-infinite sets.

\subsection{Natural density}\label{sec:dens}
Number theory offers natural density as a different notion of size that is able to make distinctions among infinite co-infinite subsets of the natural numbers. While natural density is more fine-grained than approaches based on bijection, we will review in this section that it remains coarse-grained to some extent and that the measure is not total on $\mathcal{P}(\mathbb{N})$.

The natural density (or asymptotic density) $d$ of a subset of natural numbers is defined as
\begin{equation}\label{eq:natdensity}
    d(S) \stackrel{\textrm{def}}{=} \lim_{n \rightarrow \infty} \frac{|(S \cap \{1, 2, 3, \ldots, n\})|}{n},
\end{equation}
with $S \in \mathcal{P}(\mathbb{N})$ such that this limit is indeed defined. 
In words, the natural density of a subset of $\mathbb{N}$ is the limit (if it exist at all) of the (finite) cardinality of the intersection of this set with an initial segment of $\mathbb{N}$, as this initial segment goes to infinity.

The natural density can also be understood in terms of characteristic sequences. 
The characteristic sequence of a set $S \in \mathcal{P}(\mathbb{N})$ is defined as follows, for all natural numbers $n$:\footnote{Throughout this article, we will indicate a sequence $a: \mathbb{N} \rightarrow X$ by its value at a generic position $n$, $a_n$. $\chi_n$ is a first example of this gloss.}
$$\chi_n(S) \stackrel{\textrm{def}}{=} |S \cap \{n\}|.$$
This is a binary sequence that indicates whether or not $n$ is in the set, respectively by value 1 or 0.
Now we define the sequence of partial sums of the characteristic sequence of $S$, which is a non-decreasing integer sequence:
$$f_n(S) \stackrel{\textrm{def}}{=} \sum_{i=1}^n \chi_i(S).$$

In what follows, it will be helpful to also define the non-decreasing sequence of finite intersections of a given set, $S \subset \mathbb{N}$, with initial segments of $\mathbb{N}$ of length $n$:
$$S_n \stackrel{\textrm{def}}{=} (S \cap \{1, 2, 3, \ldots, n\}).$$
The sequence of finite sizes of $S_n$ is equal to the aforementioned sequence of partial sums of the characteristic sequence: $f_n(S) = |S_n|$.

Using this notation, we may rewrite the natural density as:
$$d(S)=\lim_{n \rightarrow \infty} f_n(S)/n,$$
where $1/n$ is a normalization factor.\footnote{Observe that the natural density of a subset $S$ can be regarded as the C{\'e}saro limit of its characteristic sequence: $$C-\lim \chi_n(S) \stackrel{\textrm{def}}{=} \lim_{n \rightarrow \infty} \frac{1}{n} \sum_{i=1}^n \chi_i(S) = \lim_{n \rightarrow \infty} \frac{f_n(S)}{n} = d(S).$$\label{fn:Cesaro}}

For subsets of $\mathbb{N}$, $\lim_{n \rightarrow \infty} f_n(S)$ suffices to indicate the cardinality: it results in the finite cardinality for finite sets and it diverges for sets of cardinality $\aleph_0$.
Below, we will encounter the non-normalized sequence in other approaches as well.

The above sequences all take the canonical order of $\mathbb{N}$ for granted. Observe that, if we considered a sequence $f'_n$ that maps $n$ to the finite intersections in a different order, the resulting limit of $f'_n(S)/n$ might be different.

\subsubsection{Examples of natural densities}
The natural density is zero for all finite subsets and equal to one for all infinite co-finite subsets. Moreover, the natural density is also zero for some infinite co-infinite sets, such as the set of perfect squares, cubes, and higher powers, as well as the set of primes. So, the natural density coarse-grains finite sets and some infinite co-infinite sets.
Sets with the same natural density can be distinguished in terms of the rate of convergence. For instance, the natural density of the primes is zero, but its rate of convergence is well-studied, too: the prime number theorem implies that the natural density of primes converges to zero as $1/ln(n)$ \citep[see, e.g.,][Ch.~4]{FineRosenberger:2016}.

The complements of zero-density sets have natural density one.
The natural density does distinguish between other infinite co-infinite sets, such as the even numbers, the numbers divisible by three, four, etc., and other arithmetic progressions (i.e., translations of the former). In general, for any set
$$\mathbb{M}_{a,i} \stackrel{\textrm{def}}{=} \{n \mod a = i \hspace{0.5em}\mid\hspace{0.5em} n \in \mathbb{N} \},$$
the natural density is:
$$d(\mathbb{M}_{a,i}) = 1/a.$$
Throughout the paper, we pay special attention to the case of the set of even numbers, $\mathbb{E} \stackrel{\textrm{def}}{=} \mathbb{M}_{2,0}$, and the set of odd numbers, $\mathbb{O }\stackrel{\textrm{def}}{=} \mathbb{M}_{2,1}$, both of which have natural density 1/2.

Like the cases with a natural density of zero or one, also all intermediate assignments are coarse-grained, in the sense that there are infinitely many other sets with the same natural density. To see this, it suffices to consider the countably infinite collection of sets that make a finite difference with the given set.
Of course, taking the union or intersection of a given set with an infinite set of zero density also illustrates this, provided that the resulting set has a defined natural density. This is not guaranteed, since the sets for which $d$ is defined do not form a $\sigma$-algebra citep[pp.~23--24]{Kubilius:1964}.

\subsubsection{Subsets without a natural density\label{sec:nonatdens}}
An example of a set that does not have a well-defined natural density is the set of natural numbers whose binary expansion has an odd length, for which the `lower density' (corresponding to the lower limit in eq.~\ref{eq:natdensity}) is $1/3$ and the upper density (upper limit thereof) is $2/3$. In general, a limit point of $f_n(S)/n$ is the limit of a converging subsequence \citep[see, e.g.][]{Tao:2017}, so the lower and upper density are the extremal values of the collection of limit points of $f_n(S)/n$. To illustrate that the interval spanned by the lower and upper density need not be symmetric around $1/2$, as the previous example might suggest, consider the natural numbers that start with decimal 1, which has a lower density of $1/9$ and an upper density of $5/9$ \citep[p.~415]{Tenenbaum:2015}.

Moreover, it is possible to construct sets with lower density zero and upper density one. This can be achieved by alternating between increasingly long intervals of natural numbers that are excluded or included, such that there are subsequences of $f_n(S)/n$, indexed by the ends of either the excluded or the included intervals, that converge to zero or one, respectively.

Throughout this paper, we will focus on an example of this kind: the set $\mathbb{S}$, with super-exponential stretches. We define $\mathbb{S}$ recursively as follows: $\mathbb{S}_1 \stackrel{\textrm{def}}{=} \{\}$ and for all $n>1$, $\mathbb{S}_n \stackrel{\textrm{def}}{=} \mathbb{S}_{n-1}$ if $k \stackrel{\textrm{def}}{=} \lceil\log_2( \log_2{n}) \rceil$ is even and $\mathbb{S}_n \stackrel{\textrm{def}}{=} \mathbb{S}_{n-1} \cup \{n\}$ if $k$ is odd.\footnote{$\lceil \cdot \rceil$ indicates the ceiling function. An initial fragment of the set is given explicitly by $\{3,4,17,18,19,20\ldots\}$.}

The extreme values of $f_n(\mathbb{S})/n$ occur at $n=2^{(2^k)}$. If $k$ is odd, we have:
$$f_n^{\textrm{o}}(\mathbb{S})/n = \sum_{l=0}^k (-1)^{l+1}2^{(2^l-2^k)}.$$
Each term is smaller than the next and the signs alternate. The largest term, which occurs for $l=k$, equals $+1$; the second largest term, for $l=k-1$, is $- \frac{1}{n^{1/2}}$; and the third largest, for $l=k-2$, is $+ \frac{1}{n^{3/4}}$. Hence we see that:
$$1 - \frac{1}{n^{1/2}} < f_n^{\textrm{o}}(\mathbb{S})/n < 1 - \frac{1}{n^{1/2}} + \frac{1}{n^{3/4}}.$$
Since $\frac{1}{n^{1/2}}$ and $\frac{1}{n^{3/4}}$ tend to zero in the limit of $n$ to infinity, this shows that the limit of $f_n(\mathbb{S})/n$ tends to 1 for the subsequence where $n=2^{(2^k)}$ and $k$ is odd.\footnote{Numerical values of $f_n^{\textrm{o}}(\mathbb{S})$ at $n=2^{(2^k)}$ with $k$ odd can be found in the OEIS as the odd-numbered values of sequence A325912 \citep{Manyama:2019}.}
So, the upper density of $\mathbb{S}$ is 1.

If $n=2^{(2^k)}$ and $k$ is even, we have:
$$f_n^{\textrm{e}}(\mathbb{S})/n = \sum_{l=0}^{k-1} (-1)^{l+1}2^{(2^l-2^k)}.$$
There is one fewer term than in the odd case: this corresponds to the fact that now we are considering an initial fragment that ends on a segment of which the elements are not included. As before, each term is smaller than the next and the signs alternate. This time, the largest term, for $l=k-1$, equals $\frac{1}{n^{1/2}}$, while the second largest, for $l=k-2$, is $-\frac{1}{n^{3/4}}$. Hence, we obtain:
$$\frac{1}{n^{1/2}} - \frac{1}{n^{3/4}} < f_n^{\textrm{e}}(\mathbb{S})/n < \frac{1}{n^{1/2}}.$$
This shows that the limit of $f_n(\mathbb{S})/n$ tends to 0 for the subsequence where $n=2^{(2^k)}$ and $k$ is even. So, the lower density of $\mathbb{S}$ is 0. Since it differs from the upper density, this proves that $\mathbb{S}$ has no natural density. We will return to this set below, in the context of numerosity theory.

Natural density is ideally suited for sets with a characteristic sequence that is periodic eventually: for those, it equals the relative fraction of included elements within the period. The natural density may be undefined for other sets, even if they are recursively enumerable (i.e., Turing-recognizable), as the above examples show.

\subsubsection{Extension to $\mathbb{N}$ and relation to probability\label{sec:extendND}}
As already mentioned, the domain of $d$ does not exhaust $\mathcal{P}(\mathbb{N})$ and it is not an algebra.\footnote{\citet{KerkvlietMeester:2016a} made a proposal to extend natural density uniquely to a larger subset of $\mathcal{P}(\mathbb{N})$. While we do not consider it here, for our purposes it is important that their extension is not total.} However, it is possible to extend the natural density to a measure that is total on $\mathcal{P}(\mathbb{N})$ \citep[as discussed, e.g., in][§3.2]{WenmackersHorsten:2013}. This is achieved by replacing the standard limit in the definition (eq.~\ref{eq:natdensity}) by a generalised limit, which is a different algebra homeomorphism: a real-valued, free ultrafilter limit. For sets for which the standard natural density is undefined, this extension yields a non-unique result in the interval spanned by the lower and upper density \citep[for details, see][]{SchurzLeitgeb:2008}. Its construction crucially relies on the Hahn--Banach theorem, which requires the Boolean prime ideal theorem \citep[see, e.g.,][Ch.~12]{Schechter:1997}. So, although the full Axiom of Choice is not required, this extension does require a non-constructive axiom, which is the source of its non-uniqueness.\footnote{The Hahn--Banach theorem is sufficient to prove the Banach--Tarski decomposition of a three-dimensional sphere, which is widely regarded as pathological \citep{Pawlikowski:1991}.}
For $\mathbb{S}$, this means that its generalized limit could be 0, 1, or some intermediate number, depending on the free ultrafilter used for the limit.

Remark that the underdetermination that shows up in the generalised density for infinite co-infinite sets such as $\mathbb{S}$ does not indicate a pathology of these sets. Also, there is nothing arbitrary or underdetermined about these subsets in themselves. What is arbitrary is to choose one limit point above all others, a choice that is achieved by a non-constructive object, such that all such assignments become arbitrary but fixed at once.

Generalised density is related to the standard notion of probability \citep[without countable additivity; see, e.g.,][p.~24]{Kubilius:1964}: it is real-valued, non-negative, finitely additive, and normalized on $\mathbb{N}$. Indeed, it is an instance of the general notion of a limiting relative frequency, which has become a textbook definition of probability. As such, it can be considered as a finitely additive probability measure that is singleton-uniform on $\mathbb{N}$ (i.e., $d(F)=0$ for all finite $F \subset \mathbb{N}$).

Like all other merely finitely additive probability measures on infinite domains, generalised density is an intangible. In section~\ref{sec:num}, we will encounter more intangibles: free ultrafilters and the $\alpha$-numerosities that depend on them.

At this point, despite the explicit example of $\mathbb{S}$, perhaps it could be speculated that other sets that fail to have a natural density are random in some sense and thereby connected to intangibles. The next section shows that this speculation is misguided.

\subsubsection{Almost all subsets of $\mathbb{N}$ are intangibles with density $1/2$\label{sec:random}}
Despite the previous focus on extending the measure via a generalised limit that lacks uniqueness, almost all subsets of $\mathbb{N}$ do have a unique natural density.
To see this, first observe that almost all binary sequences are algorithmically random \citep{Martin-Lof:1966}, where `almost all' means that they have measure one in the Lebesgue measure on the Cantor space $\{0,1\}^\mathbb{N}$ (i.e., the fair Bernoulli measure on the set of infinite binary sequences).\footnote{By this canonical measure, every element of $\mathcal{P}(\mathbb{N})$ has measure zero: it treats each subset of $\mathbb{N}$ on a par (unlike $d$, for instance), regardless of its properties, so it is not an alternative measure of size.} Let us call subsets of $\mathbb{N}$ with a characteristic sequence that is algorithmically random `random subsets' \citep[cf.][§3.1]{Axon:2010}. Then, almost all subsets of $\mathbb{N}$ are random sets. Algorithmically random sequences are normal, which implies in the case of binary sequences that they have a limiting relative frequency of $1/2$. Hence, random sets have a natural density of $1/2$.
This proves our initial statement: almost all subsets of $\mathbb{N}$ have a unique natural density, namely $1/2$.

When defining randomness in terms of \textit{Kollektivs}, \citeauthor{vonMises:1928}'s (\citeyear{vonMises:1928}) guiding concept was precisely that no admissible place selection rule should pick out a subsequence with a lower or higher limiting relative frequency than 1/2. If those existed, then it would be possible to bet on random outcomes with a higher than 1/2 probability of winning, which goes against the pretheoretical concept of random outcomes. Place selection rules have later been replaced by more precise conditions on randomness, including algorithmic randomness \citep[also known as Martin-L\"of randomness;][]{Martin-Lof:1966}.

Looking at the definition of natural density (eq.~\ref{eq:natdensity}), however, it may appear strange that random subsets have a well-defined natural density. How would one compute it without an algorithmic description of the set? Indeed, an explicit computation is not possible, because random subsets are intangibles: we can demonstrate that random subsets exist in ZFC, in abundance, without being able to give a single explicit example.\footnote{This may be another reason for rejecting the notion of arbitrary subsets of $\mathbb{N}$, as \citet{Feferman:1999} did.} In this sense, the problem of determining the natural density of a fully specified random subset cannot arise. If we consider $R$ to be a random set, then $\chi_n(R)$ is not computable, so the limit in equation~\ref{eq:natdensity} for the natural density cannot be computed explicitly either, yet it must be 1/2 because $\chi_n(R)$ of a random set $R$ is normal.

\begin{figure}[!htb]
\centering
  \includegraphics[width=0.8\textwidth]{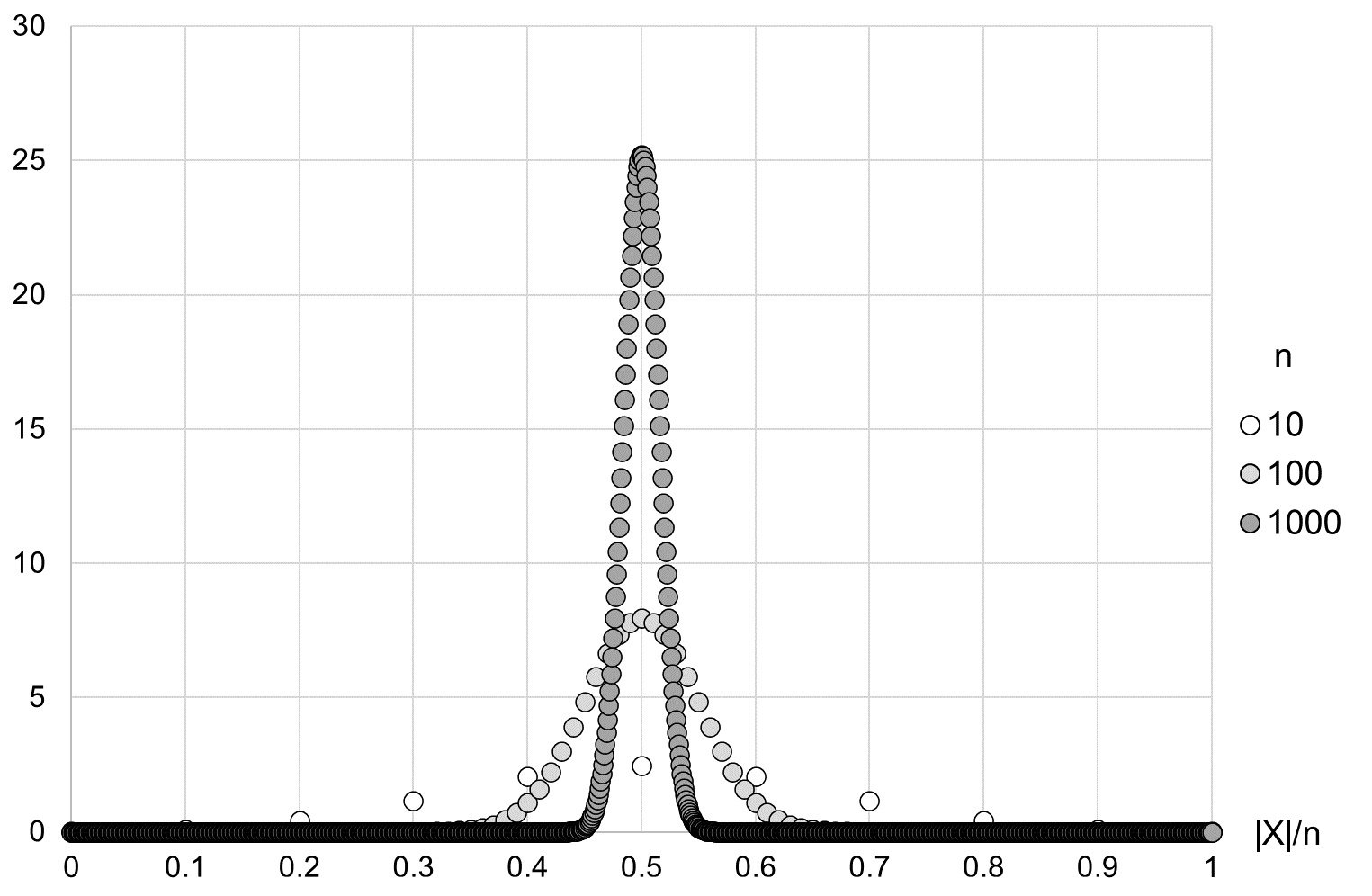}\\
  \caption{Histogram of the number of subsets of a set $F$ with a given fraction of the total number of elements, $|F|=n$, for $n$ equal to 10, 100, and 1000. This shows that the subsets, $X \subset F$, with $|X|/n=1/2$ dominate in the limit of $n$ to infinity.}\label{Fig:Binomial}
\end{figure}

Moreover, it may appear strange that natural density is measured on a continuous scale ($[0,1]$), yet almost all sets have the same value ($1/2$). This may be understood as follows.
First consider a finite set, $F$, with a number of elements, $|F|=n$, that we assume to be even for simplicity.
For each subset $X$ of $F$, we can define its $F$-density as: $d_F(X)=|X|/n$, which takes values in $V=\{0,1/n,\ldots,1/2-1/n,1/2,1/2+1/n,\ldots,1-1/n,1\}$.
Now consider the distribution that expresses how many subsets of $F$ have an $F$-density equal to each of these values, which can be visualized as a histogram: see Figure~\ref{Fig:Binomial} for some numerical examples. For each $n$, the distribution has a maximum that occurs at $d_F(X)=1/2$, corresponding to the collection of subsets with $n/2$ elements. Their number is given by the central binomial coefficient, $n!/((n/2)!)^2$, while the total number of subsets of $F$ is $2^n$ and the number of bins equals $|V|=n+1$.
So, their relative frequency evolves as $n!/((n/2)!)^2 \times (n+1)/2^n$, which diverges in the limit of $n$ to infinity. In the limit of $n$ to infinity, the full distribution goes to a delta function with peak at density 1/2 (cf.\ Figure~\ref{Fig:Binomial}).
The limiting case shows once more that almost all subsets of $\mathbb{N}$, with infinitely long characteristic sequences, have a natural density of $1/2$. These include $\mathbb{E}$ and $\mathbb{O}$, but almost all others are random sets.

\subsection{$\alpha$-Numerosity}\label{sec:num}
If we return to the initial trilemma, cardinality chose the Humean principle and natural density retained neither the Humean nor the Euclidean principle. So, we have yet to explore a theory that chooses the latter.
$\alpha$-Numerosity theory indeed takes the Euclidean principle as its starting point. $\alpha$-Numerosities of subsets of $\mathbb{N}$ (canonically ordered) can be defined axiomatically as a function $\num$ from $\mathcal{P}(\mathbb{N})$ to some set of values $\mathcal{N} \supset \mathbb{N}$ on which the total order relation $<$ and the operation $+$ are defined, as follows:\footnote{This definition is a simplified version of the more general definition of $\alpha$-numerosity on a class of labelled sets presented by \citet[Ch.~15]{BenciDiNasso:2019}; see also \citet{BenciDiNasso:2003b} for an earlier version, as well as \citet{Benci-etal:2006b} for a related approach, in which the Euclidean principle follows from the axioms.}
\begin{description}
    \item[Unit] $\forall n \in \mathbb{N}, \hspace{0.3em} \num(\{n\})=1$.
    \item[Additivity] $\forall S, T \in \mathcal{P}(\mathbb{N})$, \hspace{0.3em} if $S \cap T = \emptyset$ then $\num(S \cup T) = \num(S) + \num(T)$.
    \item[Finite approximation] $\forall S, T \in \mathcal{P}(\mathbb{N}),$ \hspace{0.3em} if $f_n(S) \leq f_n(T)$ for all $n \in \mathbb{N}$ then $\num(S) \leq \num(T)$.
    \item[Euclidean principle] $\forall S, T \in \mathcal{P}(\mathbb{N})$, \hspace{0.3em} if $S \subset T$ then $\num(S) < \num(T)$.
\end{description}

Agreement with finite cardinality follows from the first two axioms. They also ensure the Euclidean principle for finite subsets. The third axiom implies that the order on the numerosities of sets agrees with that of the finite cardinalities of their intersections with initial segments of $\mathbb{N}$; this ensures a certain harmony with natural density, as we will see below. The final axiom is needed to extend the Euclidean principle to infinite sets. In doing so, Hume's principle no longer holds, though a weaker statement applies to $\alpha$-numerosities: if $\num(S)=\num(T)$, then there exists a bijection between $S$ and $T$ (so $|S|=|T|$). The inverse direction is not guaranteed; consider, for example, $S=\mathbb{N}$ and $T$ any infinite strict subset. Then both have cardinality $\aleph_0$, yet $\num(S)>\num(T)$. So, $\alpha$-numerosity is more fine-grained than cardinality.\footnote{\citet{LynchMackey:2023} recently proposed a similar formalism called `magnum theory', which maps $\mathcal{P}(\mathbb{N})$ ``to a subclass of the surreals, the surnatural numbers {\bf Nn} (the non-negative omnific integers)''. Their theory is also designed to adhere to the Euclidean principle and one way to define the magnum is via the natural density.} Cardinality is compatible with infinitely many finer grained notions of size. Although one might have hoped that the axioms for $\alpha$-numerosity pick out one particular fine-graining, this turns out not to be the case, as we will see in section~\ref{sec:numnotunique}.
Before showing that, we first examine a model of the theory.

\subsubsection{The numerosity of $\mathbb{N}$ and the non-Archimedean $\alpha$-limit}
It follows from the axioms that $\num(\mathbb{N})$ cannot be finite and that $\num(\mathbb{N})$ is the largest numerosity value among all sets in $\mathcal{P}(\mathbb{N})$. It is customary to set
$$\num(\mathbb{N}) \stackrel{\textrm{def}}{=} \alpha;$$ hence the name $\alpha$-numerosities. 
A model for these axioms is obtained by setting the codomain $\mathcal{N}$ equal to an infinite hyperfinite set:
$$\mathcal{N}_\alpha \stackrel{\textrm{def}}{=} \{0, 1, 2, \ldots, \alpha-2, \alpha-1, \alpha\}.$$ This is an infinite initial segment of a set of hypernatural numbers\footnote{A set of hypernatural numbers, $^\ast\mathbb{N}$, indicates a non-standard model of Peano arithmetic, first discovered by \citet{Skolem:1934}.} extended with zero, $^\ast\mathbb{N} \cup \{0\}$, on which the relation $<$ (total order) and the operation $+$ are defined by transfer from the corresponding operations on $\mathbb{N}$.
The values in $\mathcal{N}_\alpha$ suffice to assign $\alpha$-numerosities to all elements of $\mathcal{P}(\mathbb{N})$. For each finite subset, the $\alpha$-numerosity is a finite number and, for each co-finite set, it is $\alpha$ minus a finite number. The $\alpha$-numerosity of an infinite co-infinite set is an infinite number that is infinitely smaller than $\alpha$. For example, for the set of even numbers, the $\alpha$-numerosity is usually taken to be $\alpha/2$ (although $\alpha/2-1$ is admissable, too, as we will discuss below).

$\alpha$-Numerosity can be understood in terms of a non-Archimedean limit operation, called the $\alpha$-limit (written as $\lim_{n \uparrow \alpha}$).
The $\alpha$-limit of a sequence can be thought of as follows: first, the sequence with indices in $\mathbb{N}$ is extended to a hypersequence with indices in $^\ast\mathbb{N}$; then, the extended sequence is evaluated at $\alpha$. Using this notation, we can write (for any $S \in \mathcal{P}(\mathbb{N}))$: $$\num(S) = \lim_{n \uparrow \alpha} f_n(S).$$
We return to this in section~\ref{sec:limits}.

\subsubsection{$\alpha$-Numerosity is totally ordered and highly non-unique\label{sec:numnotunique}}
In the preamble to their axioms, \citet{BenciDiNasso:2003b} assumed that the order on $\mathcal{N}$ is total; this indeed applies to $\langle \mathcal{N}_\alpha, < \rangle$. We will investigate the effects of taking $\mathcal{N}$ to be a merely partially ordered set in section~\ref{sec:c-num}.

So, $\alpha$-numerosities are elements of a hyperfinite set that is totally ordered. They are constrained by the partial order of the subset relation on $\mathcal{P}(\mathbb{N})$ (via the final axiom: the Euclidean principle), as well as by the other axioms, but these are not sufficient to specify a total order.
Any partial order can be extended to a total order \citep[by the Szpilrajn extension theorem; also mentioned by][]{Mancosu:2009}, but doing so requires a non-constructive choice principle (namely, the axiom of finite choice, weaker than the Axiom of Choice but also independent of ZF). For our case, this means that the axioms do not suffice to determine a total order on all of $\mathcal{P}(\mathbb{N})$ uniquely.\footnote{In the context of probability theory, a similar point has been discussed by \citet{Easwaran:2014} and \citet{Hofweber:2014}.}

To bridge the gap from a partial to a total order, the construction of $\alpha$-numerosities crucially relies on a non-constructive object: a free ultrafilter on $\mathbb{N}$, which is an intangible \citep[§6.33]{Schechter:1997}. There are uncountably many free ultrafilters,\footnote{In fact, the set of free ultrafilters on $\mathbb{N}$ has the same cardinality as $\mathcal{P}(\mathcal{P}(\mathbb{N}))$; see \citet[§6.33]{Schechter:1997} for references to proofs.} so this completion of the order is highly non-unique. In fact, using the reflections from section~\ref{sec:random}, we see that almost all members of a free ultrafilter on $\mathbb{N}$ are random subsets of $\mathbb{N}$ (which are intangibles in their own right). They have density $1/2$, as do their complements, which are also random sets. Therefore, for almost every member of a free ultrafilter, it is indeed completely arbitrary to include it rather than its complement.

The axioms of \citet{BenciDiNasso:2003b} require a specific type of free ultrafilter, called a \emph{selective} ultrafilter, the existence of which is independent of ZFC, but we do not go into that here. For our purposes, it is merely important to note that this restriction does not lead to a unique assignment.

To illustrate the non-uniqueness, let us return to the infinite co-finite set $\mathbb{S}$ that was recursively defined in section~\ref{sec:nonatdens}.
The $\alpha$-numerosity of this set depends on the ultrafilter. In particular, the most extreme values are obtained when $\alpha$ is of the form $2^{(2^\kappa)}$. Here, $\kappa$ is an infinite hypernatural, infinitely smaller than $\alpha$. There are two ways in which this can happen: either the set $\{ 2^{(2^k)} \mid k \in 2\mathbb{N}-1 \}$ is in the ultrafilter, such that $\kappa$ is odd, or $\{ 2^{(2^k)} \mid k \in 2\mathbb{N} \}$ is in the ultrafilter and $\kappa$ is even. (Observe that these sets indicate the indices of the subsequences corresponding with respectively the upper and lower density of $\mathbb{S}$.)

If $\kappa$ is odd, then $\num(\mathbb{S})$ is maximal:
$$\num(\mathbb{S}) = \lim_{n \uparrow \alpha} f_n(\mathbb{S}) = \lim_{n \uparrow \alpha} \sum_{l=0}^k (-1)^{l+1}2^{(2^l)}.$$
The result of this is a hyperfinite sum with $\alpha+1$ terms.\footnote{Hyperfinite sums generalize the standard notion of a finite sum to a sum with a hypernatural number of terms; for details see, e.g., \citet[§12.7]{Goldblatt:1998}.} It equals a particular value in $\mathcal{N}_\alpha$, but I am not aware of a closed formula to express it with.\footnote{This is not due to the $\alpha$-limit, but to the fact that a closed formula for $\sum_{l=0}^k (-1)^{l+1}2^{(2^l)}$ is unlikely to be found; in any case, \citet{Manyama:2019} reported none.} Still, we can give an indication of the value by considering the three largest terms of the sum: $\alpha$, $-\sqrt{\alpha}$, and $\sqrt{\sqrt{\alpha}}$. 
Hence,
$$\alpha-\sqrt{\alpha} < \num(\mathbb{S}) < \alpha-\sqrt{\alpha}+\sqrt{\sqrt{\alpha}}.$$
This shows that the largest possible $\alpha$-numerosity that can be assigned to $\mathbb{S}$ is smaller than that of any co-finite set, yet bigger than some other infinite co-infinite sets, such as that of non-squares (which can be chosen to have numerosity $\alpha-\sqrt{\alpha}$).

If $\alpha$ is of the form $2^{(2^\kappa)}$ with even $\kappa$, then $\num(\mathbb{S})$ is minimal:
$$\num(\mathbb{S}) = \lim_{n \uparrow \alpha} f_n(\mathbb{S}) = \lim_{n \uparrow \alpha} \sum_{l=0}^{k-1} (-1)^{l+1}2^{(2^l)}.$$
The result is a hyperfinite sum with $\alpha$ terms, the two largest of which are $\sqrt{\alpha}$ and $-\sqrt{\sqrt{\alpha}}$. Hence, in this case:
$$\sqrt{\alpha}-\sqrt{\sqrt{\alpha}} < \num(\mathbb{S}) < \sqrt{\alpha}.$$
This shows that the smallest possible $\alpha$-numerosity that can be assigned to $\mathbb{S}$ is larger than that of any finite set, yet smaller than that of some infinite co-infinite sets, such as that of the perfect squares (which can be chosen to have numerosity $\sqrt{\alpha}$).

To recap, while the $\alpha$-numerosity assignment to finite and co-finite sets is fully specified by the axioms, $\alpha$-numerosities of infinite co-infinite sets depend on the properties of the free ultrafilter used to construct them. This may result in infinite differences in the $\alpha$-numerosity of infinite co-infinite sets that do not have a natural density and finite differences for those that do (as shown in section~\ref{sec:numvsnatdens}).

\subsubsection{Critical responses to $\alpha$-numerosity\label{sec:numcritique}}
So far, we have seen that $\alpha$-numerosity theory requires an intangible that cannot be proven to exist in ZFC. Moreover, if this object exists, it is highly non-unique and particular assignments crucially depend on its properties.
These aspects play a major role in the critical reception of Euclidean theories of set size.

\citet{Godel:1947} argued that Cantor's definition of cardinality is uniquely well-motivated and the only correct notion of set size. For decades, this was the received view that nearly stopped the search for alternative conceptions of set size in its tracks.
However, \citet{Mancosu:2009} reviewed historical as well as contemporary theories on which sizes of subsets of natural numbers can be compared meaningfully, thereby contesting G{\"o}del's view that Cantor's cardinality theory was inevitable.
In particular, $\alpha$-numerosity theory was developed by \citet{BenciDiNasso:2003b}.

\citet{Parker:2013} admitted that the numerosity theory of \citet{BenciDiNasso:2003b} is logically consistent (provided that ZFC is), but he offered some new counterarguments. He argued that (total) Euclidean theories of set size are too arbitrary, because they violate weak invariance principles (such as translation or rotation invariance). This argument mainly relied on point sets in metric spaces, but for number sets (like the one under consideration here) \citet[§5]{Parker:2013} similarly concluded that they are epistemically not very useful, because not all specific size assignments are uniquely well-motivated. As a result, \citet{Parker:2013} argued that numerosity assignments are uninformative, and even misleading: since at least some of them could just as well be different, they do not reveal a stable property intrinsic to those sets. This seems to violate a deeply rooted assumption about our size conception, perhaps its very essence. So, this objection is closely linked to the non-uniqueness of the free ultrafilter, and its consequences on $\alpha$-numerosity assignments, which cannot even be fully specified because free ultrafilters are intangibles.

So, while the preamble stipulates that the order must be total, the axioms do not fully determine which infinite co-infinite sets are bigger than which. This depends on the particular selective ultrafilter.

The fact that selective ultrafilters are independent of ZFC should be flagged as a separate issue. Depending on one's philosophical views on mathematics, this may elicit different reactions, akin to those raised in response to the Continuum Hypothesis, which has a similar status.\footnote{I am grateful to a reviewer, who suggested making this analogy.} Indeed, the existence of selective ultrafilters is implied by the Continuum Hypothesis \citep[as also remarked by][§6]{Parker:2013}.
Mathematical pluralists might take the independence of ZFC to mean that it is admissible to investigate the implications of $\alpha$-numerosity theory, while it is equally admissible to explore the consequences of the negation of this possibility. Mathematical realists, on the other hand, might demand stronger, independently motivated axioms to settle the truth or falsehood of the existence of selective ultrafilters and thus of the (in)correctness of the $\alpha$-numerosity approach to set sizes.
So, it is unclear what a realist in the sense of \citet{Godel:1947} should conclude: on the one hand, they view Cantor's cardinality as the ultimate theory of set size; yet, if they complete G{\"o}del's  program and find independently justified axioms that settle the Continuum Hypothesis---and if it is found true---that implies the type of object that fuels an alternative theory of set size.

In section~\ref{sec:c-num}, we discuss a related approach by \citet{Trlifajova:2024} based on the Fr{\'e}chet filter that may escape these criticisms.

\subsubsection{Relation to natural density and additional constraints\label{sec:numvsnatdens}}
As mentioned, the third axiom of $\alpha$-numerosity theory guarantees a certain agreement with natural density. In particular, for any set $S$ such that $d(S)$ is defined, it holds that $d(S) = st(num(S)/\alpha)$, where $num(S)/\alpha$ is a hyperrational number (i.e., an element of $^\ast\mathbb{Q}$) and $st$ is the standard part function that maps a finite hyperrational number to the unique nearest standard real number \citep[for details, see][§16.9]{BenciDiNasso:2019}. This implies that sets which have the same natural density, cannot have numerosities that differ by more than a finite number.\footnote{This can be seen as follows. For any sets $S$ and $T$ such that $d(S)=d(T)$, it holds that $st(num(S)/\alpha) = st(num(T)/\alpha)$, hence $st(num(S)/num(T)) = 1$. If their ratio is infinitely close to 1, the numerosities cannot differ by more than a finite number.} 
In particular, considering random subsets, there are infinitely more subsets of $\mathbb{N}$ with an $\alpha$-numerosity within a finite interval around $\alpha/2$ than anywhere else in the hyperfinite set $\mathcal{N}_\alpha$.

In addition to the axioms, some further size comparisons can be motivated and achieved by restricting the ultrafilter accordingly. For instance, the axioms do not fix whether various sets $\mathbb{M}_{a,i}$ for given $a$ have exactly the same $\alpha$-numerosity; they could also differ by a finite amount. However, one may further stipulate \citep[cf.][pp.~63--64]{BenciDiNasso:2003b}:
\begin{align*}
\forall a \in \mathbb{N}, \forall i \in \{0, \ldots, a-1\} \hspace{1em} & \num(\mathbb{M}_{a,i}) =  \alpha/a; \\
\forall p \in \mathbb{N} \hspace{1em} & \num(\{ n^p \hspace{0.5em}\mid\hspace{0.5em} n \in \mathbb{N} \}) =  \sqrt[p]{\alpha}.
\end{align*}
We could add further stipulations to fix, for instance, the $\alpha$-numerosity of the set of natural numbers whose binary expansion has an odd length. Since this set does not have a natural density (recall that its lower and upper natural density were resp.\ $1/3$ and $2/3$), its $\alpha$-numerosity may differ by an infinite amount, depending on the ultrafilter.
Still, uncountably many $\alpha$-numerosity assignments are left to the specific properties of the (arbitrary but fixed) free ultrafilter.

We may call $num(S)/\alpha$ the $\alpha$-density of $S$.
Like natural density, $\alpha$-density can be given a probabilistic interpretation. Next, we comment on a related formalism.

\subsubsection{Relation to probability}
$\alpha$-Numerosity theory is closely related to non-Archimedean probability (NAP) theory, as discussed by \citet[§5]{Benci-etal:2013} and \citet[§3.6]{Benci-etal:2018}. The non-uniqueness of NAP functions has been discussed by \citet[§6.1]{Benci-etal:2018}. For the present purposes, we limit ourselves to the case of a NAP function that represents a singleton-uniform probability measure on $\mathbb{N}$, which is equal to a $\alpha$-numerosity assignment up to a normalization factor $1/\alpha$. So, this leads to normalized measure of sizes of all elements of $\mathcal{P}(\mathbb{N})$.

The non-uniqueness of NAP assignments has been discussed in detail for the subset of even numbers, $\mathbb{E}$, and that of odd numbers, $\mathbb{O}$, by \citet[pp.~141--142]{Benci-etal:2013}.
In terms of $\alpha$-numerosity, this boils down to two possible assignments: either $\num(\mathbb{E})=\alpha/2=\num(\mathbb{O})$, or $\num(\mathbb{E})=\alpha/2-1 < \num(\mathbb{O})=\alpha/2+1$. This difference depends on whether $\mathbb{E}$ is in the free ultrafilter used to generate the hypernatural numbers, or $\mathbb{O}$ is; or, equivalently, whether $\alpha$ is even or odd \citep[p.~63]{BenciDiNasso:2003b}. The possibility of ending up with $\num(\mathbb{E}) < \num(\mathbb{O})$, with a difference of 1, can be understood by considering that we defined $\mathbb{N}$ as starting from 1, so the odd numbers are one step ahead at the odd positions in the canonical order on $\mathbb{N}$; this constant difference remains visible in the non-Archimedean limit at $\alpha$ of $f_n(\mathbb{E})$ and $f_n(\mathbb{O})$ if $\alpha$ is odd.
Similar considerations hold for other sets of the form $\mathbb{M}_{a,i}$. We return to this at the end of section~\ref{sec:combine}.

\subsection{Composite picture}\label{sec:combine}
\begin{sidewaysfigure}[!htb]
\centering
  \includegraphics[width=1\textwidth]{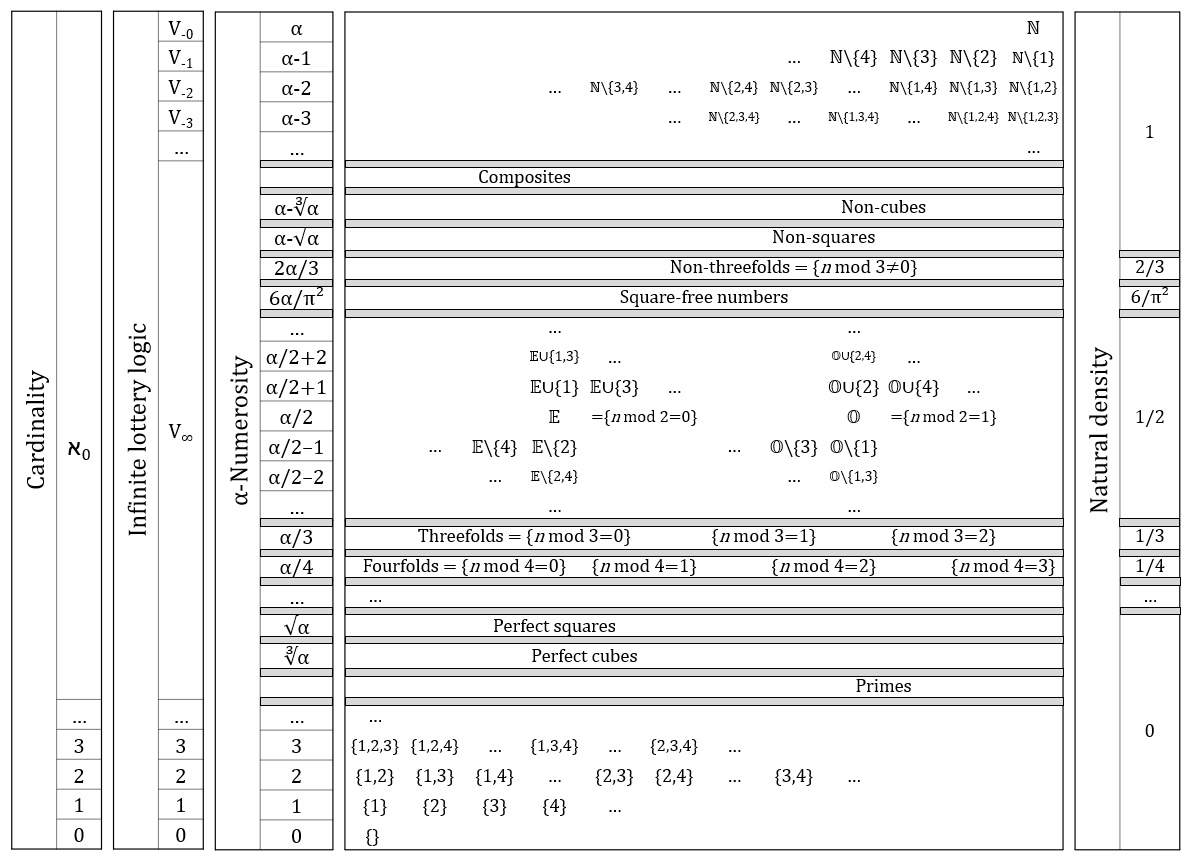}\\
  \caption{The central part of the diagram depicts part of $\mathcal{P}(\mathbb{N})$ stratified by the subset relation. Various vertical axes represent different measures of size and large cells indicate coarse-graining. Grey bands in the central diagram as well as in two scales indicate omitted parts.}\label{Fig:Assignments}
\end{sidewaysfigure}

It is time to take stock, by giving an overview of the formalisms reviewed so far.
To compare various size assignments across the formalisms and their relation to the subset lattice on $\mathcal{P}(\mathbb{N})$, they are depicted in Figure~\ref{Fig:Assignments}. Since $\alpha$-numerosities are partially underdetermined by the axioms, this picture shows just one possibility among many: one where $\alpha$ is a multiple of two, three, four, etc., as well as a power of two, three, etc. We have not included $\mathbb{S}$ in the figure, but observe that it could be placed nearly anywhere: from a point between the sets of cubes and squares until a point between the sets of non-squares and non-cubes.
We have also not included random sets: they would be at the levels of $\mathbb{E}$ and its finite difference sets.

Let us now summarize the most relevant properties of the formalisms reviewed in this section in Table~\ref{Table:comp}. Of these alternatives, $\alpha$-numerosity provides the most fine-grained measure of the sizes of subsets of natural numbers: it is total on $\mathcal{P}(\mathbb{N})$, is maximally fine-grained on finite and co-finite subsets, and not only distinguishes among infinite co-finite sets but makes similarly fine-grained distinctions among them as it does for (co-)finite ones. However, unlike the other alternatives (except for generalised density), $\alpha$-numerosity is not uniquely determined by its axioms. The Euclidean principle requires that $\alpha$-numerosity assignments respect the ordering due to the subset relation on $\langle \mathcal{P}(\mathbb{N}), \subset \rangle$. Although the third axiom restricts the possible assignments further (by demanding a form of correspondence with natural density), this still only specifies a strict partial order, while $\alpha$-numerosity values are part of the totally ordered set of hypernatural numbers, $\langle ^\ast\mathbb{N}, \subset \rangle$.
As before with generalised density, totality is provided by an intangible.

\begin{center}
\begin{table}[!htb]
\caption{Comparison of four formalisms.\label{Table:comp}}
  \begin{tabulary}{1\textwidth}{R||CCCC}
                                                 & Cardinality           & Infinite lottery logic & $\alpha$-Numerosity           & Natural density \\
\hline
\hline
Total on $\mathcal{P}(\mathbb{N})$                 & Yes                   & Yes                    & Yes                  & No, but can be extended \\
Unique assignments & Yes                   & Yes                    & No                   & Yes, but not the extension \\
Normalized (size$(\mathbb{N})=1$)                     & No, not normalizable & No, not normalizable  & No, but normalizable & Yes \\
\hline
Distinguishes between:                           &                       &                        &                      & \\
\hline
Finite sets                                        & Yes                   & Yes                    & Yes                  & No \\
Co-finite sets                             & No                    & Yes                    & Yes                  & No \\
Infinite co-infinite sets                           & No                    & No                     & Yes                  & Yes \\
  \end{tabulary}
  \end{table}
\end{center}

As reviewed in section~\ref{sec:numcritique}, in response to the results of $\alpha$-theory, \citet{Parker:2013} concluded that (total) Euclidean assignments of set size contain arbitrary aspects that are therefore misleading. He did see some hope for a different Euclidean theory of size limited to partial assignments \citep[e.g., ``{\bf Even} might be regarded as neither smaller than {\bf Odd}, nor larger, nor equal''][p.~609]{Parker:2013}, but we will see in section~\ref{sec:c-num} that partially ordered numerosity can do better than this (e.g., it excludes the possibility that $\mathbb{E}$ might be larger than $\mathbb{O}$ in a non-arbitrary way).

A similar tension between totality and arbitrariness has been discussed in the context of probability theory by \citet{Easwaran:2014}. As summarized by \citet[p.~543]{Benci-etal:2018}: ``Easwaran observes that real-valued probability functions leave out part of the structure of the partial order, whereas hyperreal-valued probability functions add structure in an arbitrary way.'' Likewise, when assigning sizes to subsets of $\mathbb{N}$, it seems that we are forced to choose between options with too little or too much structure: whereas cardinality, infinite lottery logic, and natural density do not preserve the partial order of the subset lattice, $\alpha$-numerosity comes with excess structure added by the specifics of some free ultrafilter.

\citet{Easwaran:2014} suggested using standard probability functions (which miss some structure on the partial order of the events) together with the subset structure of the algebra of events: in the case of equal probabilities, one may still treat the larger set as more probable (e.g., when deciding betting preferences). He argued that this is better than accepting the arbitrary excess structure baked into a non-standard probability function.
If we apply Easwaran's proposal to our current case, we should reject $\alpha$-numerosity in favour of natural density, while keeping an eye on the subset relation when ranking set sizes.
For example, $d(\mathbb{E} \setminus \{2\}) = d(\mathbb{E})$ but $\mathbb{E} \setminus \{2\} \subset \mathbb{E}$, so we may treat the former as `smaller', in some sense. However, this proposal does not seem sufficient, since $d(\{1\})=d(\{2, 3\})$, but $\{1\}$ and $\{2, 3\}$ are not comparable by the subset relation, so they are still `equal', in the same sense of size.

An alternative approach has been suggested, but not developed, for non-Archimedean probabilities; \citet[p.~544]{Benci-etal:2018}: ``one can always consider the entire family of NAP functions modelling a given situation, rather than an arbitrary representative of it (see also Wenmackers and Horsten [2013]). [\ldots] As a whole, the family shows us how much the probabilities of a given event, and the order of probabilities of multiple events, can vary (dependent on the choice of ultrafilter).'' This proposal is related to imprecise or interval-based probabilities. Likewise, we could consider the collection of $\alpha$-numerosity functions ranging over all selective ultrafilters.

In the next section, we review a very recent approach that manages to trace the contours of such a collection, without listing its members.

\section{Partially ordered c-numerosity of subsets of $\mathbb{N}$\label{sec:c-num}}
It is clear that all $\alpha$-numerosity assignments agree on finite and co-finite subsets of $\mathbb{N}$, and that they differ amongst each other by at least 1 on infinite co-infinite sets. Now, one may wonder about the collection of subsets that all free ultrafilters have in common. This is the co-finite filter or Fr{\'e}chet filter on $\mathbb{N}$:
$$\mathcal{F} \stackrel{\textrm{def}}{=} \{ S \subseteq \mathbb{N} \hspace{0.5em}\mid\hspace{0.5em} \mathbb{N} \setminus S \textrm{\ is \ finite} \}.$$
So, $\mathcal{F}$ consists of all co-finite sets of $\mathbb{N}$; it is free but not an ultrafilter itself \citep[see, e.g.,][p.~73]{Jech:2003}. It can be defined in ZF, so it does not require any non-constructive choice principle.
This suggests the possibility of developing a fully (quasi-)constructive theory akin to $\alpha$-numerosity.\footnote{The potential of developing non-standard analysis using the constructive co-finite filter instead of free ultrafilters has also been considered, for instance, by \citet{Palmgren:1998}.}

Recently, \citet{Trlifajova:2024} has made a reconstruction of Bolzano's work and indeed arrived at such a formalism based on the co-finite filter rather than a free ultrafilter. Like $\alpha$-numerosity theory, \citet{Trlifajova:2024}'s theory applies to all countable sets, but we limit our review to its results for subsets of $\mathbb{N}$. Again, elements of such sets are only considered in their canonical order.

\citet{Trlifajova:2024} reconstructed the work of Bolzano as being about sequences of the form $f_n(S)$ (i.e., $|S_n|$; $\sigma(S)$ in her notation), where $S$ is a subset of natural numbers. First, she showed that the set $\{f_n(S) \hspace{0.5em}\mid\hspace{0.5em} S \in \mathcal{P}(\mathbb{N})\}$ exhausts the set of non-decreasing sequences of natural numbers, which we call $\mathcal{S}$. Then, she equipped this set with addition (and multiplication, which we do not consider here) defined component-wise, and with equality and ordering defined as equality or ordering (resp.) of terms eventually.

The notion of equality or ordering of terms `eventually' means the relation holds between the terms for a co-finite subset of indices of the sequences; \textit{i.e.}, such that the set of indices is an element of the co-finite filter $\mathcal{F}$. This can be viewed as a type of limit operation, along the co-finite filter. We will call it the \textit{$\mathcal{F}$-limit} ($\lim_\mathcal{F}$).

With natural density, different sequences $f_n(S)/n$ may have the same limit: these sequences form an equivalence class under the standard limit operation. Likewise, what really matters for size attributions in the sense of c-numerosity is not an individual size sequence function, $f_n(S)$, either, but rather its equivalence class under the Fr{\'e}chet filter. Hence, it is convenient to quotient the Fr{\'e}chet filter out of the space of non-decreasing sequences: $\mathcal{S} / \mathcal{F}$.\footnote{See \citet[p.~112; emphasis added]{Trlifajova:2024}: ``The codomain of $\sigma$ is the set of nondecreasing sequences of natural numbers \emph{modulo the Fr{\'e}chet filter} which is just partially and not linearly ordered.''} With addition, equality and order defined as before, \citet[p.96]{Trlifajova:2024} showed $\mathcal{S} / \mathcal{F}$ to be a ``partially ordered non-Archimedean commutative semiring''.
This means that it contains infinite elements (as before) and retains the same arithmetical properties as $\langle \mathbb{N}, +, \times \rangle$ (which is a totally ordered Archimedean commutative semiring).

Now, equivalence classes of $f_n$ are elements of $\mathcal{S} / \mathcal{F}$. They play a role very similar to the $\alpha$-numerosity function, $\num$, but with the co-finite filter, rather than a free (and selective) ultrafilter. Therefore, we will call this the \textit{c-numerosity}, where the `c' stands for co-finite as well as for constructive:
$$\forall S \in \mathcal{P}(\mathbb{N}), \hspace{1em} \cnum(S) \stackrel{\textrm{def}}{=} \lim_\mathcal{F} f_n(S) \stackrel{\textrm{def}}{=} [f(S)]_\mathcal{F} \in \mathcal{S} / \mathcal{F}.$$

Another precursor to this approach is found in \citet{Peano:1910}:\footnote{I am grateful to a referee for providing me with a copy of this source.} he proposed to associate with each sequence an `end' value (`{\it fine}' in Italian), which can be understood as a type of limit different from the standard one. In general, \citet[p.~780]{Peano:1910} took two sequences to have the same end value if they are equal from a certain index onward. For a sequence that is constant eventually, its end value is defined to be equal to that constant. (Hence, all real values are among the ends.) For sequences such as $a_n=1/n$, however, Peano argued that their ends are actual infinitesimals and, for sequences such as $a_n=n$, actual infinities. He defined the sum and product of end values as the end value of the sum or product of the corresponding sequences. \citet{Peano:1910} did not use his construction to discuss set sizes, but applying it to sequences of the form $f_n(S)$ yields results equivalent to those of \citet{Trlifajova:2024}.\footnote{To my knowledge, \citeauthor{Peano:1910}'s (\citeyear{Peano:1910}) paper has not been translated. \citet[pp.~477--478]{BottazziKatz:2021} described the gist of it, rephrased in contemporary terms, as the construction of a ``partially ordered non-Archimedean ring'' consisting of the set of ends of real-valued sequences ``with zero divisors that extend $\mathbb{R}$''.}

Observe that c-numerosities \emph{also} provide a model for the simplified axioms for $\alpha$-numerosities that we presented in section~\ref{sec:num}, provided that we drop the requirement from the preamble that $\mathcal{N}$ is totally ordered. This amounts to setting the codomain of the new numerosity function to $\mathcal{N}_c \stackrel{\textrm{def}}{=} \mathcal{S} / \mathcal{F}$. In section~section~\ref{sec:num}, we used a totally ordered hyperfinite set as the codomain, $\mathcal{N}_\alpha = \{0, 1, 2, \ldots, \alpha-2, \alpha-1, \alpha\}$, which results in the core difference with merely partially ordered c-numerosities covered here.

C-numerosities agree with $\alpha$-numerosities on finite sets. Choosing
$$\cnum(\mathbb{N}) \stackrel{\textrm{def}}{=} \alpha,$$
the assignments also agree on co-finite sets. 
In other words, if we restrict $\mathcal{N}_c$ to the values assigned to finite and cofinite sets, that subset is totally ordered.
So, what about infinite co-infinite sets? In general, their c-numerosity (in $\mathcal{N}_c$) cannot be mapped to a unique value in $\mathcal{N}_\alpha$. Returning to the example of even and odd numbers, \citet[p.~96]{Trlifajova:2024} wrote: ``While there are one fewer or as many even numbers as odd numbers, their sum is equal to $\alpha$.'' Adapted to our notation, she obtained:
$$\cnum(\mathbb{E}) \leq \cnum(\mathbb{O}) \leq \cnum(\mathbb{E}) + 1;$$
$$\cnum(\mathbb{O}) + \cnum(\mathbb{E}) = \alpha.$$
Although $\mathbb{E}$ and $\mathbb{O}$ are comparable to each other on this approach, neither of these sets is comparable to the infinite co-infinite set $\mathbb{S}$ that was recursively defined in section~\ref{sec:nonatdens}. This example shows that c-numerosities are only partially ordered. Before we elaborate on this example, it is instructive to examine the limit operations underlying the various approaches.

\subsection{Comparing the limit operations\label{sec:limits}}
Any limit operation on sequences requires an equivalence relation that consists of two aspects: a notion of tolerance and a notion of qualified index sets, given by a filter.\footnote{For a helpful overview, see \citet{Tao:2017}. For a comparison that focuses on the difference between the standard limit and the $\alpha$-limit, see \citet[§6.1]{Wenmackers:2019a}.}
For the standard limit operation on sequences, $\lim_{n \rightarrow \infty}$, the qualified index sets are those with a sufficiently large index and the tolerance is given by arbitrarily small differences along those indices (e.g., $|a_n-a_{n-1}|<1/n$ `eventually', i.e., for all $n>N$ for some $N$). In this case, the qualified index sets are given by the co-finite sets in $\mathbb{N}$. The family of such sets forms a filter: the Fr{\'e}chet filter, $\mathcal{F}$, on $\mathbb{N}$.

A non-Archimedean limit operation, such as the $\alpha$-limit, differs from the standard limit on both accounts: the qualified index sets are given by a different type of filter, a free ultrafilter on $\mathbb{N}$ (on which additional conditions may be imposed), and the tolerance is given by exact equality (zero tolerance).

This presentation suggests two additional types of limit operations: 
one that combines a free ultrafilter with arbitrarily small differences, which was the generalised ultrafilter limit we used to extend the natural density to $\mathcal{P}(\mathbb{N})$ in section~\ref{sec:extendND}, and another one that combines the Fr{\'e}chet filter with zero tolerance, which was $\lim_\mathcal{F}$ introduced in section~\ref{sec:c-num}.

Hence, natural density, generalised density, $\alpha$-numerosity, and c-numerosity can all be understood as a type of limit operation applied to $f_n(S)$ or $f_n(S)/n$: Table~\ref{Table:limits} gives an overview.\footnote{Recall from section~\ref{sec:dens} that the standard limit of $f_n(S)$ \emph{is} defined for all sets $S$, even when that of $f_n(S)/n$ is not.}

\begin{table}[!htb]
\caption{Different formalisms in terms of the type of limit operation (rows) and sequence (columns).\label{Table:limits}}
\begin{tabular}{llr||cc}
tolerance & filter &                & $f_n$                                           & $f_n/n$                 \\
\hline
\hline
\multirow{2}{*}{$< \epsilon(n)$} & $\mathcal{F}$        & standard limit & finite cardinality & natural density         \\
                                 & $\mathcal{U}_\alpha$ & generalised limit   &   or $\infty$                               & generalised density \\
\multirow{2}{*}{$=0$}  & $\mathcal{F}$        & $\mathcal{F}$-limit        & c-numerosity                                    & c-density               \\
                                 & $\mathcal{U}_\alpha$ & $\alpha$-limit & $\alpha$-numerosity                             & $\alpha$-density       
\end{tabular}
\end{table}

We could also introduce mirror cardinalities by considering $f_n(\mathbb{N} \setminus S)$ as well. Moreover, observe that finite and co-finite sets give rise to sequences that are \textit{extremely} well-behaved in the standard limit. As remarked in footnote~\ref{fn:Cesaro}, the natural density of a set can be regarded as the C{\'e}saro limit of its characteristic function. But for any finite set, $F$, $\lim_{n \rightarrow \infty} \chi_n(F) = 0$ and $\lim_{n \rightarrow \infty} \chi_n(\mathbb{N} \setminus F) = 1$, so the additional averaging provided by the C{\'e}saro limit is not needed to get convergence in this case.

\begin{figure}[!htb]
\centering
  \includegraphics[width=0.5\textwidth]{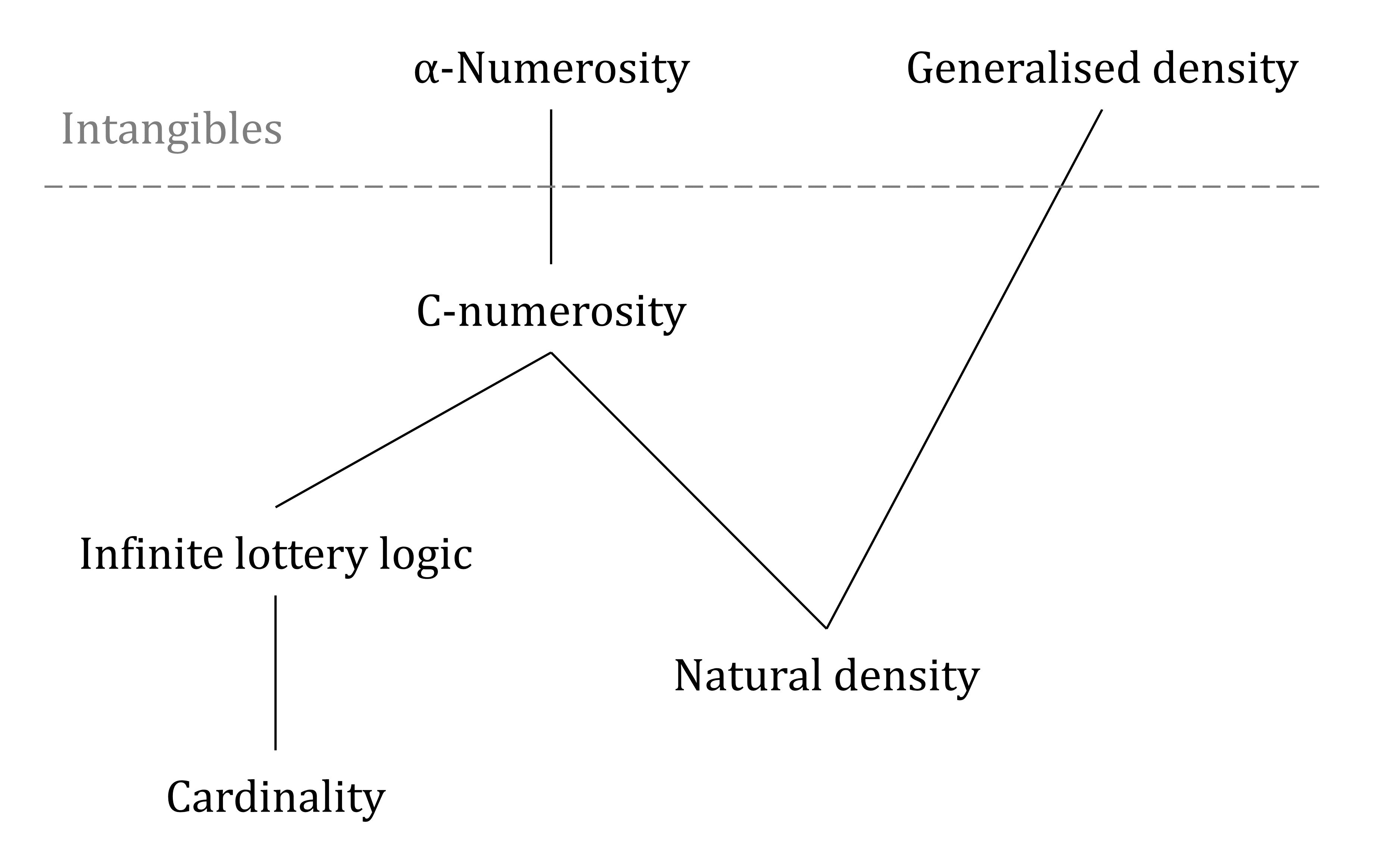}\\
  \caption{Partial order of the distinctions in set sizes made by various formalisms: more fine-grained towards the top. The dashed line indicates that intangibles occur above it.}\label{Fig:Distinctions}
\end{figure}

Table~\ref{Table:limits} helps us to understand the relations of fine-grainedness between the measures, as well as the source of partial underdetermination present in some of them.
The $\mathcal{F}$-limit of sequences is more fine-grained than the standard limit, because it requires exact equality of terms rather than equality up to an arbitrarily small difference (while using the same notion of eventuality); and it is more coarse-grained than the $\alpha$-limit because it uses a smaller filter (while using the same notion of equivalence).
In Figure~\ref{Fig:Distinctions}, we depicted the partial order on the formalisms induced by their fine-grainedness.

Both the standard limit and the $\mathcal{F}$-limit are constructive, the other two are not. Hence, natural density and c-numerosity are free of arbitrary structure introduced by intangibles, whereas both generalised density and $\alpha$-numerosity show signs of them.

\subsection{Another look at set $\mathbb{S}$}
Now that we have a better grasp of the connection between the different formalisms via the underlying limit operations, let us illustrate it with the example the set $\mathbb{S}$ (as defined in section~\ref{sec:nonatdens}).

While the standard limit of $f_n(\mathbb{S})$ is well-defined ($+\infty$), that of $f_n(\mathbb{S})/n$ is not. In other words, the natural density of $\mathbb{S}$ is undefined. As we discussed in section~\ref{sec:nonatdens}, the reason is that $f_n(\mathbb{S})/n$ has different subsequences, each with a different limit (ranging from 0 to 1). The natural density is expressed by real numbers, which are totally ordered, and the measure does not depend on arbitrary structure, but this is only possible by excluding sets like $\mathbb{S}$ from its domain.

A generalised density does assign a limit-value to $f_n(\mathbb{S})/n$, equal to that of some arbitrary but fixed subsequence. So, this measure is total on $\mathcal{P}(\mathbb{N})$ and assigns values from a totally ordered range, but with a partially arbitrary character.

Likewise, as we have seen in section~\ref{sec:numnotunique}, the $\alpha$-numerosity of $\mathbb{S}$, defined as the $\alpha$-limit of $f_n(\mathbb{S})$, can be any of a wide range of infinite numbers (depending on which subsequence is selected by the underlying ultrafilter). The fact that all admissible $\alpha$-numerosities are infinite corresponds with the divergence of the standard limit. The fact that they differ by more than a finite amount reflects the lack of a unique natural density.
In particular, the $\alpha$-limit of $f_n(\mathbb{S})/n$ varies from an infinitesimal less than $1/\sqrt{\alpha}$ (with standard part 0) to more than $1-1/\sqrt{\alpha}$, which is infinitesimally close to 1 (i.e., has standard part 1). These two assignments correspond to the two subsequences we considered before to argue that $f_n(\mathbb{S})/n$ has subsequences with limits of 0 and 1.
So, $\alpha$-numerosity is similar to general density in terms of totality, arbitrary, and lack of uniqueness. Moreover, its range is non-Archimedean to allow for the Euclidean principle.

Finally, the c-numerosity of $\mathbb{S}$ is comparable to some sets (for instance, it is larger than any finite set and also larger than the set of cubes), but not to others (such as the set of squares, $\mathbb{E}$, $\mathbb{O}$---or any of the $\mathbb{M}_{a,i}$, for that matter).
This measure is also non-Archimedean, to allow for the Euclidean principle. This time, the measure is unique and total on $\mathcal{P}(\mathbb{N})$ but now the order is partial, again to allow for sets like $\mathbb{S}$.

Recall that, if one accepts infinite co-infinite sets as unproblematic in general, there is nothing intrinsically pathological or unknowable about $\mathbb{S}$. It is not an intangible itself; its characteristic function just keeps alternating over stretches that grow at a super-exponential rate. As a result, its growth rate ($f_n(\mathbb{S})$) does not stabilize relative to that of initial fragments of $\mathbb{N}$ (of length $n$).

Most of the formalisms that we reviewed are dealing with this fact by leaving something undetermined: either by not assigning a value to some sets at all or by not totally ordering the values (natural density and c-numerosity, respectively), or by assigning them totally ordered values in some globally consistent but arbitrary way (general density and $\alpha$-numerosity).
In all these cases, something has to go: totality (whether totality on the domain or total ordering of the co-domain) or uniqueness of assignment, respectively.
The cardinality-based approaches do not face this dilemma, but they also leave unexpressed much non-arbitrary structure that the other theories agree upon regarding the sizes of infinite co-finite sets.

\subsection{The relative strength of c-numerosity}
Since the c-numerosity approach does not invoke intangibles, it escapes the criticism that \citet{Parker:2013} directed at total Euclidean theories (reviewed in section~\ref{sec:numcritique}).
However, \citet[§9]{Parker:2013} also anticipated Euclidean theories with ``partial size assignments''. While total theories were too arbitrary on his view, he expected partial theories to be too weak and narrow: unable to decide the relative sizes of some simple sets. In particular, \citet[p.~609]{Parker:2013} expected that such theories would not make any comparison between the size of the sets of even and odd numbers, $\mathbb{E}$ and $\mathbb{O}$, ranking them as ``neither smaller [\ldots], nor larger, nor equal''. However, we have seen that c-numerosity is a partial Euclidean theory that does exclude that $\mathbb{O}$ could be larger than $\mathbb{E}$. It constrains the order as follows: $\cnum(\mathbb{E}) \leq \cnum(\mathbb{O}) \leq \cnum(\mathbb{E}) + 1$ with $\cnum(\mathbb{E}) + \cnum(\mathbb{O}) = \alpha$. So, c-numerosity is much less weak than anticipated by the criticism of \citet[§9]{Parker:2013}. \citep[See also the response in][§6.]{Trlifajova:2024}

The axioms of numerosity theory merely specify a partial order. 
For numerosity values to be totally ordered, this has to be postulated separately, as \citet{BenciDiNasso:2003b} did. Doing so comes at the cost of introducing arbitrary and unknowable structure due to an intangible. C-numerosity theory, in contrast, honours the Euclidean principle and the other axioms of numerosity theory, while remaining fully constructive. It is as fine-grained as possible and as coarse-grained as needed. So, c-numerosity seems to hit the right balance between expressiveness and leaving out comparisons that are not stable.
This leads to a merely partially ordered notion of size, which means that its assignments are further removed from the usual notion of number. Moreover, the inequalities of c-numerosity assignments exactly trace the hull of all possible $\alpha$-numerosity assignments. We explore this departure from our familiar size concept and the relation between the two theories in the next section.

\subsection{Epistemicism or supervaluation\label{sec:superval}}
Taking our inspiration from formal methods for dealing with vagueness, we may consider a particular $\alpha$-numerosity function (generated by a specific free ultrafilter) as a `precisification' of c-numerosity.\footnote{Both $\alpha$- and c-numerosity quotient a filter out of the space of non-decreasing sequences: the Fréchet filter used in c-numerosity is contained in any free ultrafilter used in $\alpha$-numerosity, so the former assignments are compatible with the former, but not vice versa.}
Likewise, we may consider each generalised density as a precisifation of the natural density.
We could also say that a total extension of c-numerosity (and natural density) is multiply realizable.

In analogy to epistemicists about vagueness \citep[such as][]{Williamson:1994}, we might assume that there is one true Euclidean notion of size---given by one correct precisification, but which is unknowable to us. In this case, it is not the Euclidean notion of size itself but our knowledge of it that is incomplete. In principle, such a position should be acceptable to a mathematical realist.
On the one hand, the fact that $\alpha$-numerosity functions---like the free ultrafilters that generate them---are intangible is congenial to this interpretation: the unknowability is guaranteed. Yet, it seems hard to defend that any specific one could be privileged, since these objects do indeed represent much arbitrary structure.

Alternatively, we may view the Euclidean notion of size embodied in the axioms reviewed at the start of section~\ref{sec:num} as partially underdetermined. Then, like supervaluationists, we may treat all precisifications on a par and consider the set of such precisifications. In the study of vagueness, the notions of `supertrue' and `superfalse' are defined as truth values that agree on the set of all precisifications. Analogously, here, we may consider subsets of $\mathbb{N}$ that have a `super-$\alpha$-numerosity' as those that have the same $\alpha$-numerosity on each precisification:\footnote{Equivalently, we may say that it is supertrue that a certain set has a particular $\alpha$-numerosity. This suggestion has been made  in the context of NAP functions, by \citet[p.~544]{Benci-etal:2018}.} these are exactly the ones that have a unique c-numerosity, i.e., just the finite and co-finite sets.

One complication, mentioned before, is that we may \emph{define} the $\alpha$-numerosity function to assign $1/a$ to each subset $\mathbb{M}_{a,0}$ etc.; so, we have to evaluate the issue relative to such potential additional constraints at hand.
Observe that, although \citet{Trlifajova:2024} did not suggest this (and, in fact, argued against doing so), additional constraints may be added to partially ordered numerosities, too. (After all, adding, e.g., $\mathbb{E}$ to $\mathcal{F}$ and extending it such that it forms a filter, does not turn it into a free ultrafilter.)

All infinite co-infinite sets fail to have a super-$\alpha$-numerosity.
We may wish to distinguish cases where the $\alpha$-numerosities of different precisifications only differ by a finite amount from those where they differ by a larger margin; the former are the sets with a natural density. Mutatis mutandis, we may also consider a dual approach with subvaluations: a sub-$\alpha$-numerosity is an $\alpha$-numerosity that is assigned by at least one precisification.

Admittedly, the terms `underdetermination', `precisification', and `subvaluation' all take total measures with totally ordered values as the standard. This is indeed a standard in the sense that earlier theories that assign sizes to subsets of $\mathbb{N}$ adhere to it. Moreover, according to \citet{Forti:2022}, general comparability is an essential property of the classical notion of magnitude or size. 

\citet{Forti:2022} motivated the assumption that Euclidean set sizes are totally ordered as follows: ``Following the ancient praxis of comparing magnitudes of homogeneous objects, a very general notion of size of sets, whose essential property is general comparability of sizes, can be given through a total preordering $\preceq$ of sets according to their sizes [\ldots]'' Here, `homogeneous' means that the objects are of the same kind, with the same dimensions, though this condition was relaxed by later mathematicians. For instance, the cardinality of $\mathbb{N}^2$ can be compared to that of $\mathbb{N}$ (and they are equal). Although this concerns a preorder (or quasiorder: i.e., a binary relation that is reflexive and transitive), rather than an order (which is additionally required to be antisymmetric), it is a \emph{total} preorder, on which any two elements are comparable \citep[as][specified in his footnote 3]{Forti:2022}.

Nevertheless, once this classical assumption has been made explicit, it becomes possible to question and perhaps abandon it. A possible lesson that we could draw from studying numerosity theories is exactly that Euclidean set sizes are merely partially ordered. 
Our review shows that the notion of size has already undergone drastic changes; this may be the next, comparatively conservative, extension of it. Reluctance to accept the possibility of a merely partially ordered notion of set size may also resemble the epistemicist approach to vagueness: by positing crisp thresholds, it replaces vagueness by unknowable and non-constructive elements \citep[cf.][§3]{Sorensen:2022}. In both cases, this trade is possible, but by no means necessary.

\subsection{A final look at the size of $\mathbb{E}$ versus $\mathbb{O}$}
It may be instructive to return to our running example of $\mathbb{E}$ versus $\mathbb{O}$ one final time. Let us consider their even and numbered subsequences of $f_n/n$:
\begin{align*}
f_{n \in \mathbb{E}}(\mathbb{E}) & =1/2; \hspace{1em} & f_{n \in \mathbb{E}}(\mathbb{O}) & =1/2;\\
f_{n \in \mathbb{O}}(\mathbb{E}) & =1/2-1/(2n); \hspace{1em} & f_{n \in \mathbb{O}}(\mathbb{O}) & =1/2+1/(2n).
\end{align*}
These sequences have the same standard limit point, $1/2$; hence, the natural density of $\mathbb{E}$ and $\mathbb{O}$ is equal.
But $f_{n \in \mathbb{O}}(\mathbb{E})<f_{n \in \mathbb{O}}(\mathbb{O})$ for all $n$, so $f_n(\mathbb{E})/n$ has a strictly smaller `$\alpha$-limit point' than $f_n(\mathbb{O})/n$ if and only if $\mathbb{O}$ is in the ultrafilter. Hence, the $\alpha$-density of $\mathbb{E}$ may be smaller than that of $\mathbb{O}$, but since both $\mathbb{E}$ and $\mathbb{O}$ have natural density $1/2$, it seems arbitrary which one is in the ultrafilter to decide this matter. We might as well choose neither, as in the Fr{\'e}chet filter, which yields the inequalities of the c-density.
Since the natural density for $\mathbb{E}$ and $\mathbb{O}$ is defined, the generalised limit must equal this value, too.

Although the generalised limit involves a free ultrafilter, and thus the result may in general depend on whether or not $\mathbb{O}$ is included, the difference is glossed over in this case due to the tolerance built into this limit operation. Only for sets $S$ such that $f_n(S)/n$ has multiple standard limit points, which thus lack a natural density, the tolerance built into the generalised limit is not enough to absorb these differences.

Using the terminology of section~\ref{sec:superval}, we may say that $1/2$ and $1/2(1-1/(2\alpha-1))$ are both sub-$\alpha$-numerosities of $\mathbb{E}$ while the set lacks a super-$\alpha$-numerosity.
Put differently, it is subtrue that the set of even numbers has the same numerosity as the set of odd numbers, but this equality is not supertrue.

We have to admit that debating the $\alpha$-numerosity values of infinite co-infinite sets
may skirt dangerously close to discussing how many angels can dance on the head of a pin. \citet{Mancosu:2009} quoted Descartes (\textit{Principes de la Philosophie}, I.26) on the issue at hand: ``We will not bother to reply to those who ask if the infinite number is even or odd or similar things since it is only those who deem that their mind is infinite who seem to have to tackle such difficulties.''
Using sophisticated tools, it seems that the answer `it depends' can be refined to `it depends on which ultrafilter you use to investigate the matter,' but these answers may reveal more about the tools that were used (intangibles from $\mathcal{P}((\mathcal{P}(\mathbb{N}))$) than about the objects under study (elements of $\mathcal{P}(\mathbb{N})$) themselves. In contrast, c-numerosity refuses in a Cartesian way to answer whether the size of $\mathbb{N}$ is even or odd, or any other questions that go beyond the restrictions implied by its axioms.

\section{Conclusion and outlook}\label{sec:concl}
We have reviewed six formalisms for assigning sizes to subsets of $\mathbb{N}$: cardinality, infinite lottery logic with mirror cardinalities, natural density, generalised density, $\alpha$-numerosity and c-numerosity.
By evaluating these formalisms together and by studying the connections between them, we gained a better understanding of the essential trade-offs between different approaches to assigning sizes to subsets of $\mathbb{N}$.

We paid special attention to $\alpha$-numerosity and c-numerosity, which both respect the Euclidean principle, while providing a measure that is total on $\mathbb{N}$ and fine-graining the natural density.
$\alpha$-Numerosity is the most fine-grained of the two and yields a total order of its values. However, this comes at the cost of arbitrary structure associated with intangibles.
C-numerosity does not rely on intangibles, and hence does not introduce arbitrary structure in its notion of size. As such, it overcomes earlier criticisms raised against Euclidean theories of size \citep[in particular from][]{Parker:2013}. However, it does require us to drop another long-held assumption about set sizes, i.e., that they are totally ordered. This demand comes on top of dropping both the Humean principle and the assumption that sizes are Archimedean quantities, which all Euclidean theories of size applicable to infinite sets have to give up.

In other words, if we start from the classical dilemma between the Humean and the Euclidean principle \citep[previously discussed by][]{Mancosu:2009} and walk the Euclidean path, we encounter another fork down the road. The additional dilemma asks us to choose between accepting the consequences of invoking an intangible ($\alpha$-numerosity) or dropping the total order of sizes (c-numerosity).
Rather than arguing for either of the options, I suggest that $\alpha$-numerosity and c-numerosity are best viewed as two sides of the same coin: the inequalities of c-numerosity trace the intervals in which all possible $\alpha$-numerosity values lie. Together, they show us that there is an inherent limitation on how closely a Euclidean notion of size, applicable to all subsets of $\mathbb{N}$, can resemble its counterpart for finite sets.

Density-based approaches, which were intended to measure relative size, escape the first dilemma and show a third path (retaining neither principle). Still, they similarly face the second dilemma: a choice between invoking an intangible (generalised density) or dropping totality (natural density).

Moreover, we reconstructed the six theories as different types of limit operations on sequences of 
partial sums of the characteristic sequence of a given set, $f_n(S)$, or the corresponding density, $f_n(S)/n$. As a result, we could trace the differences in fine-grainedness of these measures as well as the non-constructive aspects present in two of them to the underlying limit operation.

Although the overview might suggest that all four combinations of two properties (considering two types of filters on $\mathbb{N}$ and two tolerance principles) have now been fully developed, there is room for exploring new variants. After all, natural density and c-numerosity are both based on the Fr{\'e}chet filter with a different tolerance principle (recall Table~\ref{Table:limits}). The fact that natural density is an Archimedean measure, while c-numerosity is not, is a direct consequence of the stricter tolerance principle of the latter theory. Yet, they \emph{also} differ in the way they give up totality: while natural density gives up totality on the domain, c-numerosity merely gives up totality of the ordering on the co-domain.

This suggests two new approaches, which have not been developed (as far as I know). First, we could consider a `gappy' numerosity theory, that only assigns values to subsets of $\mathbb{N}$ that have a super-$\alpha$-numerosity (i.e., only to finite and co-finite sets). This would be a constructive theory, weaker than c-numerosity. In fact, it would be similar to infinite lottery logic, with $V_\infty$ removed by a big gap.
Second, and more interestingly, we could consider a constructive extension of natural density, which is total on $\mathcal{P}(\mathbb{N})$, by giving up the total order on its values. Hence, this new measure cannot take real values, since they form a totally ordered field. Instead, it takes values on a merely partially ordered field. However, this means that this field cannot be complete either \citep{DeMarr:1967}, so we would have to give up the least upper bound property or one of the field axioms. To the best of my knowledge, this avenue has not been explored as an alternative for measuring set size.

Finally, there is room for Euclidean theories that do not aim for correspondence with natural density. Dropping the third axiom gives rise to a total, non-Archimedean, partially ordered measure, which does not compare infinite co-infinite sets unless they are related by the subset relation. It is stronger than infinite lottery logic but weaker than c-numerosity.

\section*{Acknowledgments}
I am grateful to Tomasz Placek and Klaas Landsman for their feedback on an earlier version and to two anonymous reviewers for their constructive reports, which helped me to improve this paper.


\bibliographystyle{plainnat}
\bibliography{Bib}

\begin{thebibliography}{45}
\providecommand{\natexlab}[1]{#1}
\providecommand{\url}[1]{\texttt{#1}}
\expandafter\ifx\csname urlstyle\endcsname\relax
  \providecommand{\doi}[1]{doi: #1}\else
  \providecommand{\doi}{doi: \begingroup \urlstyle{rm}\Url}\fi

\bibitem[Axon(2010)]{Axon:2010}
L.~M. Axon.
\newblock \emph{Algorithmically Random Closed Sets and Probability}.
\newblock PhD thesis, University of Notre Dame, Graduate Program in Mathematics, Notre Dame, Indiana, 2010.

\bibitem[Benci and {Di Nasso}(2003)]{BenciDiNasso:2003b}
V.~Benci and M.~{Di Nasso}.
\newblock Numerosities of labelled sets: {A} new way of counting.
\newblock \emph{Advances in Mathematics}, 173:\penalty0 50--67, 2003.

\bibitem[Benci and {Di Nasso}(2019)]{BenciDiNasso:2019}
V.~Benci and M.~{Di Nasso}.
\newblock \emph{Alpha-Theory: {M}athematics with Infinite and Infinitesimal Numbers}.
\newblock World Scientific, Singapore, 2019.

\bibitem[Benci et~al.(2006)Benci, {Di Nasso}, and Forti]{Benci-etal:2006b}
V.~Benci, M.~{Di Nasso}, and M.~Forti.
\newblock An {A}ristotelian notion of size.
\newblock \emph{Annals of Pure and Applied Logic}, 143:\penalty0 43--53, 2006.

\bibitem[Benci et~al.(2013)Benci, Horsten, and Wenmackers]{Benci-etal:2013}
V.~Benci, L.~Horsten, and S.~Wenmackers.
\newblock Non-{A}rchimedean probability.
\newblock \emph{Milan Journal of Mathematics}, 81:\penalty0 121--151, 2013.
\newblock \doi{10.1007/s00032-012-0191-x}.

\bibitem[Benci et~al.(2018)Benci, Horsten, and Wenmackers]{Benci-etal:2018}
V.~Benci, L.~Horsten, and S.~Wenmackers.
\newblock Infinitesimal probabilities.
\newblock \emph{British Journal for the Philosophy of Science}, 69:\penalty0 509--552, 2018.
\newblock \doi{10.1093/bjps/axw013}.

\bibitem[Bottazzi and Katz(2021)]{BottazziKatz:2021}
E.~Bottazzi and M.~G. Katz.
\newblock Infinitesimals via {C}auchy sequences: {R}efining the classical equivalence.
\newblock \emph{Open Mathematics}, 19:\penalty0 477--482, 2021.

\bibitem[Cantor(1895)]{Cantor:1895}
G.~Cantor.
\newblock Beitr{\"a}ge zur {B}egr{\"u}ndung der transfiniten {M}engenlehre.
\newblock \emph{Mathematische Annalen}, 46:\penalty0 481--512, 1895.

\bibitem[Dedekind(1888)]{Dedekind:1888}
R.~Dedekind.
\newblock \emph{Was sind und was sollen die {Z}ahlen?}
\newblock Vieweg, Braunschweig, Germany, 1888.

\bibitem[DeMarr(1967)]{DeMarr:1967}
R.~DeMarr.
\newblock Partially ordered fields.
\newblock \emph{The American Mathematical Monthly}, 74:\penalty0 418--420, 1967.

\bibitem[Easwaran(2014)]{Easwaran:2014}
K.~Easwaran.
\newblock Regularity and hyperreal credences.
\newblock \emph{Philosophical Review}, 123:\penalty0 1--41, 2014.

\bibitem[Feferman(1999)]{Feferman:1999}
S.~Feferman.
\newblock Does mathematics need new axioms?
\newblock \emph{The American Mathematical Monthly}, 106:\penalty0 99--111, 1999.

\bibitem[Fine and Rosenberger(2016)]{FineRosenberger:2016}
B.~Fine and G.~Rosenberger.
\newblock \emph{Number Theory. An Introduction via the Density of Primes}.
\newblock Birkh{\"a}user, Cham, Switzerland, 2016.
\newblock Second edition.

\bibitem[Forti(2022)]{Forti:2022}
M.~Forti.
\newblock A {E}uclidean comparison theory for the size of sets.
\newblock \url{http://arxiv.org/abs/arXiv:2212.05527}, 2022.

\bibitem[Frege(1884)]{Frege:1884}
G.~Frege.
\newblock \emph{Die Grundlagen der Arithmetik}.
\newblock Verlag von Wilhelm Koebner, Breslau, 1884.

\bibitem[Galilei(1638)]{Galileo:1638}
G.~Galilei.
\newblock \emph{Discorsi e dimostrazioni matematiche, intorno {\`a} due nuove scienze}.
\newblock Elsevier, Leiden, 1638.
\newblock Translated by H.\ Crew and A.\ de Salvio, introduced by A. Favaro ``Dialogues Concerning Two New Sciences''. Macmillan New York, 1914.

\bibitem[G{\"o}del(1947)]{Godel:1947}
K.~G{\"o}del.
\newblock What is {C}antor's continuum problem?
\newblock \emph{American Mathematical Monthly}, 54:\penalty0 515--525, 1947.

\bibitem[Goldblatt(1998)]{Goldblatt:1998}
R.~Goldblatt.
\newblock \emph{Lectures on the Hyperreals; An Introduction to Nonstandard Analysis}, volume 188 of \emph{Graduate Texts in Mathematics}.
\newblock Springer, New York, NY, 1998.

\bibitem[Hofweber(2014)]{Hofweber:2014}
T.~Hofweber.
\newblock Infinitesimal chances.
\newblock \emph{Philosophers' Imprint}, 14:\penalty0 1--14, 2014.

\bibitem[Hume(1739--1740)]{Hume:1739}
D.~Hume.
\newblock \emph{A Treatise of Human Nature}.
\newblock London, UK, 1739--1740.
\newblock Reprinted by Oxford University Press, Oxford, UK, 2000.

\bibitem[Jech(2003)]{Jech:2003}
T.~J. Jech.
\newblock \emph{Set Theory}.
\newblock Springer Monographs in Mathematics. Springer, Berlin, Germany, 2003.
\newblock Reprint of the Third Millennium Edition, Revised and Expanded.

\bibitem[Kerkvliet and Meester(2016)]{KerkvlietMeester:2016a}
T.~Kerkvliet and R.~Meester.
\newblock Uniquely determined uniform probability on the natural numbers.
\newblock \emph{Journal of Theoretical Probability}, 29:\penalty0 797--825, 2016.

\bibitem[Kubilius(1964)]{Kubilius:1964}
J.~Kubilius.
\newblock \emph{Probabilistic Methods in the Theory of Numbers}, volume~11 of \emph{Translations of Mathematical Monographs}.
\newblock American Mathematical Society, Providence, Rhode Island, 1964.
\newblock Translated by G.\ Burgie and S.\ Schuur.

\bibitem[Lynch and Mackey(2023)]{LynchMackey:2023}
P.~Lynch and M.~Mackey.
\newblock Counting sets with surreals. {P}art {I}: {S}ets of natural numbers.
\newblock arXiv preprint arXiv:2311.09951, 2023.

\bibitem[Mancosu(2009)]{Mancosu:2009}
P.~Mancosu.
\newblock Measuring the size of infinite collections of natural numbers: {W}as {C}antor's theory of infinite number inevitable?
\newblock \emph{The Review of Symbolic Logic}, 2:\penalty0 612--646, 2009.

\bibitem[Mancosu(2015)]{Mancosu:2015}
P.~Mancosu.
\newblock In good company? {O}n {H}ume's principle and the assignment of numbers to infinite concepts.
\newblock \emph{The Review of Symbolic Logic}, 8:\penalty0 370--410, 2015.

\bibitem[Manyama(2019)]{Manyama:2019}
S.~Manyama, 2019.
\newblock \url{https://oeis.org/A325912}; retrieved May 6, 2024.

\bibitem[Martin-L{\"o}f(1966)]{Martin-Lof:1966}
P.~Martin-L{\"o}f.
\newblock The definition of random sequences.
\newblock \emph{Information and Control}, 9:\penalty0 602--619, 1966.

\bibitem[Maschio(2020)]{Maschio:2020}
S.~Maschio.
\newblock Natural density and probability, constructively.
\newblock \emph{Reports in Mathematical Logic}, 55:\penalty0 41--59, 2020.

\bibitem[Norton(2021)]{Norton:2021}
J.~D. Norton.
\newblock Eternal inflation: {W}hen probabilities fail.
\newblock \emph{Synthese}, 198:\penalty0 3853--3875, 2021.

\bibitem[Palmgren(1998)]{Palmgren:1998}
E.~Palmgren.
\newblock Developments in constructive nonstandard analysis.
\newblock \emph{The Bulletin of Symbolic Logic}, 4:\penalty0 233--272, 1998.

\bibitem[Parker(2013)]{Parker:2013}
M.~W. Parker.
\newblock Set size and the part--whole principle.
\newblock \emph{Review of Symbolic Logic}, 6:\penalty0 589--612, 2013.

\bibitem[Pawlikowski(1991)]{Pawlikowski:1991}
J.~Pawlikowski.
\newblock The {H}ahn--{B}anach theorem implies the {B}anach--{T}arski paradox.
\newblock \emph{Fundamenta Mathematicae}, 138:\penalty0 21--22, 1991.

\bibitem[Peano(1910)]{Peano:1910}
G.~Peano.
\newblock Sugli ordini degli infiniti.
\newblock \emph{Rendiconti della Reale Accademia dei Lincei}, 19:\penalty0 778--781, 1910.
\newblock Reprinted in: G.\ Peano, \textit{Opere Scelte}, vol.\ 1, Edizioni Cremonese, Roma (1957), pp.\ 359--362.

\bibitem[Schechter(1997)]{Schechter:1997}
E.~Schechter.
\newblock \emph{Handbook of Analysis and its Foundations}.
\newblock Academic Press, San Diego, CA, 1997.

\bibitem[Schurz and Leitgeb(2008)]{SchurzLeitgeb:2008}
G.~Schurz and H.~Leitgeb.
\newblock Finitistic and frequentistic approximation of probability measures with or without $\sigma$-additivity.
\newblock \emph{Studia Logica}, 89:\penalty0 257--283, 2008.

\bibitem[Skolem(1934)]{Skolem:1934}
Th.~A. Skolem.
\newblock {\"U}ber die {N}icht-charakterisierbarkeit der {Z}ahlenreihe mittels endlich oder abz{\"a}hlbar unendlich vieler {A}ussagen mit ausschliesslich {Z}ahlenvariablen.
\newblock \emph{Fundamenta Mathematicae}, 23:\penalty0 150--161, 1934.

\bibitem[Sorensen(2022)]{Sorensen:2022}
R.~Sorensen.
\newblock Vagueness.
\newblock In E.~N. Zalta and U.~Nodelman, editors, \emph{Stanford Encyclopedia of Philosophy}. 2022.
\newblock \url{https://plato.stanford.edu/archives/win2023/entries/vagueness/}.

\bibitem[Tao(2017)]{Tao:2017}
T.~Tao.
\newblock Generalisations of the limit functional, 2017.
\newblock URL \url{https://terrytao.wordpress.com/2017/05/11/generalisations-of-the-limit-functional/}.

\bibitem[Tenenbaum(2015)]{Tenenbaum:2015}
G.~Tenenbaum.
\newblock \emph{Introduction to Analytic and Probabilistic Number Theory}, volume 163 of \emph{Graduate Studies in Mathematics}.
\newblock American Mathematical Society, 2015.

\bibitem[Trlifajov{\'a}(2024)]{Trlifajova:2024}
K.~Trlifajov{\'a}.
\newblock Sizes of countable sets.
\newblock \emph{Philosophia Mathematica}, 32:\penalty0 82--114, 2024.

\bibitem[von Mises(1928)]{vonMises:1928}
R.~von Mises.
\newblock \emph{Probability, Statistics and Truth}.
\newblock George Allen \& Unwin, London, UK, 1928.
\newblock Reprinted by Dover, Mineola, NY, 1981.

\bibitem[Wenmackers(2019)]{Wenmackers:2019a}
S.~Wenmackers.
\newblock Infinitesimal probabilities.
\newblock In J.~Weisberg and R.~Pettigrew, editors, \emph{Open Handbook of Formal Epistemology}, pages 199--265. PhilPapers Foundation, 2019.
\newblock URL \url{https://philpapers.org/archive/WENIP.pdf}.

\bibitem[Wenmackers and Horsten(2013)]{WenmackersHorsten:2013}
S.~Wenmackers and L.~Horsten.
\newblock Fair infinite lotteries.
\newblock \emph{Synthese}, 190:\penalty0 37--61, 2013.

\bibitem[Williamson(1994)]{Williamson:1994}
T.~Williamson.
\newblock \emph{Vagueness}.
\newblock Routledge, London, UK, 1994.

\end{thebibliography}

\end{document}